\documentclass[12pt]{article}

\usepackage{amssymb,amsbsy,amsmath}

\usepackage{graphicx} 
\usepackage{Topsy-style-arxiv}  

\begin{document}

\newcommand{\TableAlg}{
\begin{table}[t]
\OBSCURE{
\Dm[width=5.5em,height=2em]
\Agoth         & \Fcal\Agoth& \Mcal\Agoth             & \Fcalbar\Agoth \\
\JBB{A},\JDD{A}&\JBB{\Sblock},\JDD{\Sblock}
                            &\JBB{\alpha},\JDD{\alpha}&             \\
A              & \Sblock    &  \OS                    & \Sblockbar  \\
\dTo~\alpha    &\dTo~{\proj}&\dTo~{\id}               &\dTo~{\proj} \\
\OS            &  S         & \OS                     & S           \\
\mD
}
\begin{center}
\includegraphics[scale=0.67,trim=0pt 0pt 0pt 0pt,clip]{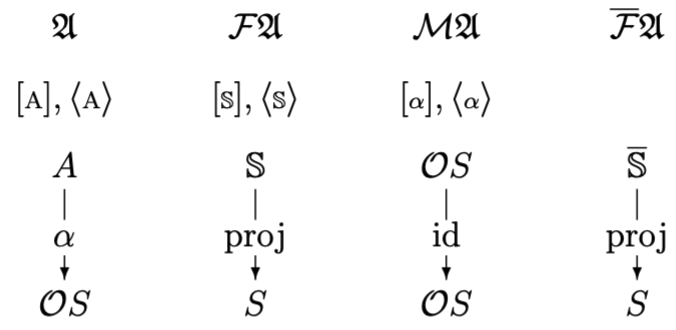}
\end{center}
\caption{The furnishing for a garden and its manicured version}
\label{table-alg}
\centerline{\rule{0.95\textwidth}{0.20ex}}\vspace{2ex}
\end{table}}

\newcommand{\TableGeo}{
\begin{table}[t]
\OBSCURE{
\Dm[width=5.5em,height=2em]
\Sgoth         & \Gcal\Sgoth&\Hcal\Sgoth &\Gcalbar\Sgoth \\
\JBB{\Sblock},\JDD{\Sblock}
               &\JBB{\sigma},\JDD{\sigma}
                            &\JBB{\Sblockbar},\JDD{\Sblockbar}
                                         &\JBB{\sigmaB},\JDD{\sigmaB} \\
\Sblock        & \OS        &\Sblockbar  & \OS          \\
\dTo~\sigma    &\dTo~{\id}  &\dTo~{\proj}&\dTo~{\id}    \\
S              & \OS        & S          & \OS          \\
\mD
}
\begin{center}
\includegraphics[scale=0.67,trim=0pt 0pt 0pt 0pt,clip]{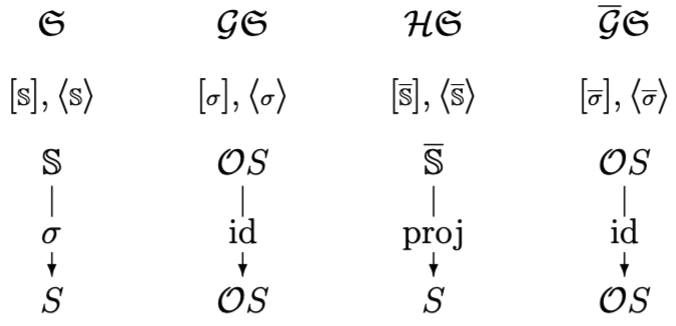}
\end{center}
\caption{The furnishing for a plot and its harvested version}
\label{table-geo}
\centerline{\rule{0.95\textwidth}{0.20ex}}\vspace{2ex}
\end{table}}

\newcommand{\TableNat}{
\begin{table}[t]
\OBSCURE{
\Dm[height=1.5em,width=5ex]
\Sblock &         &\rTo^{\GunitS}&         &\Fblock &\hspace{4ex}&
P       &\rMapsto&\big(\sigma(P),\St(P),\Bl(P)\big)               \\
        &\rdTo~{\sigma}&       &\ldTo~{\proj} &        &&
\dMapsto&        &\dMapsto               \\
        &              &S      &              &        \\
\dTo^{\Phi}&       &\dTo~{\phi}&&\dTo~{(\phi,\,\phi_\ast,\,\phi^\La)}\\
        &              &T      &              &        \\
        &\ruTo~{\tau}  &       & \luTo~{\proj}&        &&
        &  &
\big((\phi\circ\sigma)(P),(\phi_\ast\circ\St)(P),(\phi^\La\circ\Bl)(P)\big) \\
\Tblock &        &\rTo_{\Gunit_\Tgoth}&          &\Gblock &&
\Phi(P) &\rMapsto&
\big((\tau\circ\Phi)(P),  (\St\circ\Phi)(P), (\Bl\circ\Phi)(P)\big)   \\
\mD
}
\begin{center}
\includegraphics[scale=0.67,trim=0pt 0pt 0pt 0pt,clip]{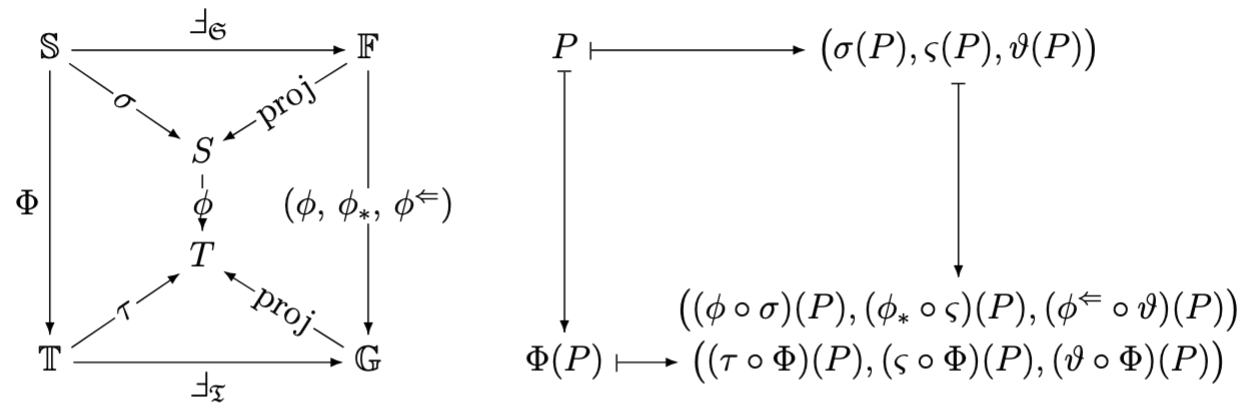}
\end{center}
\caption{A naturality square}
\label{table-205}
\centerline{\rule{0.95\textwidth}{0.20ex}}\vspace{2ex}
\end{table}}

\title{Topologically valued transition structures}
\author{
Matthew Collinson\\ 
\small{School of Natural and Computing Sciences},
University of Aberdeen,\\
\small{Aberdeen, AB24 3UE, United Kingdom,}
\texttt{\small{matthew.collinson@abdn.ac.uk}}
}
\date{}
\maketitle

\begin{abstract}
We investigate several categories related to transition structures,  
using a mixture of algebraic and topological methods. We
show how two such categories are connected by a contravariant
adjunction. This is the most detailed of a family 
of such results depending on topological restrictions on 
objects and morphisms. 
\end{abstract}

\bigskip

A transition structure consists of a set of nodes furnished with a 
binary relation called the transition relation. 
Such structures, with a variety of names, arise in diverse parts of 
mathematics, from dynamical systems to logic. More recently, 
they have become one of the fundamental tools of theoretical computer 
science. 
 
Two extreme forms of morphism between transition structures are commonly
considered. 
On one hand there are the obvious relation preserving 
functions. These are often too weak to be of use. 
On the other hand 
there are those relation preserving functions 
which are also simulations. This is a very strong condition which is 
described in the next section. 
There have been several attempts to find morphisms
intermediate between these extremes using topolological 
restrictions. More often than not, the topology is 
imposed directly on the set of nodes, following the tradition of 
\cite{hal55}, see for example \cite{esakia}, \cite{gold76a}, 
\cite{samvac}, \cite{Wij90}.  
A more recent example is \cite{hil00}, which can be seen as a special 
case of the present work. 

In this paper we depart from that tradition. We attach a 
topological space to each transition structure. The nodes of the
structure and the points of the space may be distinct, but there is a
function which assigns a point to each node, thought of as the value
of the node. Thus each path through the transition relation produces a
sequence of points in the space. 
The separation of nodes from points gives great flexibility. For
instance, we can investigate a single transition structure valued over
different spaces, or we may study a whole family of transition
structures valued over a single space. The same kind of structure has recently been studied under the name `empirical variable' \cite{bvb}, but with a different intended logical application in mind.

We first set up these topological transition structures and related 
algebraic gadgets. This leads to several categories, but we concentrate 
on just two of these, $\PlotL$
and $\Grdn$, which seem to be the most important. We investigate the
properties of these categories, and conclude by showing they are
contravariantly adjoined.

{\textbf{Acknowledgement: }} This paper was produced through joint work with the late Harold Simmons. It is identical to a draft produced circa 2006 except for the requirement of surjectivity of valuation function in the definition of topological transition structure (plot), the remark that surjectivity is used in the proof of Theorem~\ref{item-220-j}, and a citation \cite{bvb} with accompanying remark. 

\section{Preliminaries}\label{sec1}

We first review the standard notions and results on which we
build. We then introduce the structures and morphisms we consider
throughout the paper. 

\Defn\label{item-210-a}
A {\sf transition structure} is a set $\Sblock$ furnished
with a binary relation $\rra$, the transition relation.

A {\sf transition morphism}
\OBSCURE{
\Dm \Sblock&\rTo^\Phi&\Tblock \mD
}
\begin{center}
\includegraphics[scale=0.73]{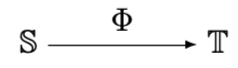}
\end{center}
between a pair of transition structures is a function
$\Phi:\Sblock\rTo\Tblock$ such that 
$$P\rra Q \Rra \Phi(P) \rra \Phi(Q)$$
holds for all nodes $P,Q\in\Sblock$. 

A {\sf simulation} from $\Sblock$ to $\Tblock$ is a transition morphism $\Phi$
such that for each node $P\in\Sblock$ and transition $\Phi(P)\rra R$ in
$\Tblock$, there is a transition $P\rra Q$ in $\Sblock$ with
$\Phi(Q)=R$ \ED 
\EDefn

These structures are known by various names. In particular, they are
sometimes called `Kripke frames', which is a bit unfortunate for a
`frame' is also a quite different structure that we need
here. These other frames are discussed later starting
with Definition \ref{new1.5}.

We refer to a member $P\in\Sblock$ as a {\sf node}. 
This enables us to reserve `point' and
`element' for other things. The transition relation $\rra$ on
$\Sblock$ is always hidden. Thus we refer to a `transition
structure' $\Sblock$. Furthermore, when there are two transition
structures around, as in Definition \ref{item-210-a}, we use $\rra$
for both the carried transition relations. 

It is routine to check that transition structures and transition
morphisms form the objects and the arrows of a category. However, 
these morphisms are too weak for most uses of
transition structures. A simulation is nothing more than a functional
bisimulation. Simulations together with transition structures form a 
category, but in this case the arrows are very strong. Notions
intermediate between that of morphism and simulation are needed. 
In this paper we show how to obtain such notions. 

\Defn
Let $\Sblock$ be a transition structure

For each node $P\in\Sblock$ we set
$$P^\ra = \{Q\in\Sblock\,|\,P\rra Q\}$$
to obtain the set of nodes $P$ can reach in one step.

For each $E\subseteq\Sblock$ we use
$$P\in\JBB{\Sblock}(E) \LRa P^\ra\subseteq E \hspace{5ex}
  P\in\JDD{\Sblock}(E) \LRa P^\ra\mbox{ meets }E$$
for $P\in\Sblock$ to obtain two operators 
$\JBB{\Sblock}$ and $\JDD{\Sblock}$ on the power set $\Pcal\Sblock$.   \ED
\EDefn

Each of these operators is the dual complement of the other, that is
$$\JDD{\Sblock}(E)=(\JBB{\Sblock}(E'))'$$ 
for each $E\subseteq\Sblock$. Each operator
determines the parent transition relation, since we have
$$P\rra Q \LRa P\in\JDD{\Sblock}(\{Q\})$$
for each $P,Q\in\Sblock$.
Thus the study of a particular transition structure
$\Sblock$ is equivalent to the study of one or both of the operators
$\JBB{\Sblock}, \JDD{\Sblock}$ on $\Pcal\Sblock$. Of course, such
operators have special properties, as shown by the following well known
characterization.

\Lem
A companion pair $\JBB{\Sblock},\JDD{\Sblock}$ of operators on a power set
$\Pcal\Sblock$ arise from a transition relation on $\Sblock$ if and
only if 
$$\JBB{\Sblock}\left(\bigcap\Ecal\right) = 
\bigcap\left\{\JBB{\Sblock}(E)\,|\,E\in\Ecal\right\}
\hspace{5ex}
  \JDD{\Sblock}\left(\bigcup\Ecal\right) = 
\bigcup\left\{\JDD{\Sblock}(E)\,|\,E\in\Ecal\right\}$$
hold for each family $\Ecal$ of subsets of $\Sblock$.
\ELem

An obvious idea is to turn these infinitary notions into finitary
ones. One could replace the power set $\Pcal\Sblock$ by a boolean
algebra $A$ and allow this to carry a dual complementary pair of operators
$\JBB{A},\JDD{A}$ where the first is a $\wedge$-semilattice morphism
and the second is a $\vee$-semilattice morphism. We follow a
different path. We first replace the power set
$\Pcal\Sblock$ by a topology (thus moving from the discrete to the
continuous), and then replace the topology by a `frame' in the other
sense (thus moving from the point-sensitive to the point-free).

We assume known the idea of a topological space $S$ with its topology
$\OS$ of open sets. Two spaces $S$ and $T$ are compared
by a continuous map
\OBSCURE{
\Dm S&\rTo^\phi& T \mD
}
\begin{center}
\includegraphics[scale=0.73]{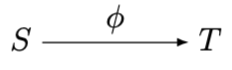}
\end{center}
and so we get the category of spaces and maps.
Every function has an associated inverse image function (going in the
opposite direction). The characteristic property of a continuous map
$\phi$, as above, is that it produces an inverse image function
\OBSCURE{
\Dm \OS&\lTo^{\phi^\la}&\OT \mD
}
\begin{center}
\includegraphics[scale=0.73]{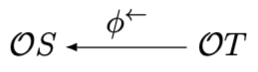}
\end{center}
between the carried topologies. This has certain algebraic properties
(which are discussed in the next subsection).
We write
$$(\cdot)\tclose  \hspace{5ex} (\cdot)\topen \hspace{5ex} (\cdot)'$$
for the closure operation, the interior operation, and the complement
operation carried by the power set $\Pcal S$ of a space $S$. There
are also other bits of gadgetry that we need.

\Defn
Let $S$ be a space.

The {\sf specialization order} on $S$ is the relation $\specord$ given
by
$$p\specord q \LRa p\in q\tclose$$
for points $p,q\in S$.

For each $E\subseteq S$ we use
$$q\in E\tsat \LRa (\exists p\in E)[p\specord q]$$
(for $p\in S$) to obtain the {\sf saturation} $E\tsat$ of $E$.

For each $E\subseteq S$ we use
$$E\tlens = E\tsat\cap E\tclose$$
to obtain the {\sf lens closure} $E\tlens$ of $E$.
\EDefn

The specialization order is a pre-order. It is a partial order
precisely when the space is $T_0$, and is equality precisely when the
space is $T_1$.
The saturation of a set is the upper section in the
specialization order generated by the set.
A lens of a space is a subset $E$ with $E\tlens=E$. These are certain
convex parts of the space in the specialization order. 

Point-sensitive topology is routine and doesn't need much
explanation. 
Point-free topology is less well known. Fortunately, we don't need
many particular details, and all of these can be described quite
quickly. 

The use of the word `frame' here must not be confused with the use in
`Kripke frame'. They are very different notions. 

\Defn\label{new1.5}
A {\sf frame} $(A,\leq,\bigvee,\wedge,\bot,\top)$ is a complete lattice
which satisfies the frame distributive law
$$a\wedge\bigvee X = \bigvee\{a\wedge x\,|\,x\in X\}$$
for all $a\in A$ and $X\subseteq A$. A frame morphism (between frame)
is a function that preserves the distinguished attributes, those being 
arbitrary suprema and finitary infima, including the top and
bottom. \ED
\EDefn

Each frame morphism $f^\ast=f:B\rTo A$ has a right adjoint $f_\ast$
\OBSCURE{
\Dm B&\pile{\rTo^{f^\ast} \\ \\ \lTo_{f_\ast}}& A \mD
}
\begin{center}
\includegraphics[scale=0.73,trim=0pt 0pt 0pt 0pt,clip]{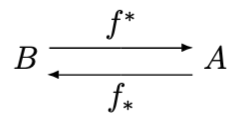}
\end{center}
satisfying
$$f^\ast(b) \leq a \LRa b\leq f_\ast(a)$$
for $b\in B$ and $a\in A$. This right adjoint is a useful technical
device. It is a monotone function but
need not be a frame morphism. 
We sometimes write $f^\ast$ for
the frame morphism $f$ when we make use of the right adjoint.

For every space $S$ the topology $\OS$ is a frame (with the evident
operations). For each continuous map $\phi:S\rTo T$ the inverse image
function $\phi^\la:\OT\rTo\OS$ is a frame morphism. This gives  a
contravariant functor $\mathcal{O}$ from the category of spaces to the
category of frames. This contravariance is the reason for our
apparently perverse notation.

Frames are a way of getting at the
algebraic properties of a topology $\OS$ without referring to the
points of the parent space $S$. In this paper we also
capture some transition properties in an algebraic way. The following
are the appropriate notions. 

\Defn
Let $A$ be a frame.

An {\sf operator} on $A$ is a monotone function $A\rTo A$.

A {\sf $\Box$-operator} on $A$ is an operator $\JBB{A}$ which is a
$\wedge$-semilattice morphism, that is
$$\JBB{A}(\top)=\top \hspace{5ex} 
\JBB{A}(x\wedge y) = \JBB{A}(x) \wedge \JBB{A}(y)$$
for all $x,y\in A$.

Given a $\Box$-operator on $A$, a {\sf companion $\Diamond$-operator}
is an operator $\JDD{A}$ on $A$ such that
$$\JBB{A}(x) \wedge \JDD{A}(y)  \leq \JDD{A}(x\wedge y)$$
for all $x,y\in A$.

A {\sf bed} is a frame furnished with a companion pair of operators. \ED
\EDefn

In Subsection
\ref{subsec2.3} we generate many substantial examples of
enriched frames
so for now we make do with a few rather pathological ones.

\Egs
Let $A$ be any frame.

(a) There are two extreme $\Box$-operators given by
$$\JBB{A}(x)= x \hspace{5ex} 
  \JBB{A}(x)=\left\{\begin{array}{@{}c@{\mbox{ if }}c}
           \top&x=\top \\[0.5ex] \bot&x\neq\top
           \end{array}\right.$$
(for $x\in A$).  

(b) The operator given by
$$\JDD{A}(x)=\bot$$
(for $x\in A$) is a companion $\Diamond$-operator for every
$\Box$-operator carried by $A$.

(c) The identity operator given by
$$\JDD{A}(x)=x$$
(for $x\in A$) is a companion $\Diamond$-operator for every
deflationary $\Box$-operator carried by $A$.
\EEgs

These examples show that a $\Box$-operator can have more than one
companion $\Diamond$-operator, and each $\Diamond$-operator can be the
companion of more than one $\Box$-operator. This is true even when the
frame is boolean.

Beds are compared using frame morphisms
which interact with the furnishings in an appropriate way. 
However, the form of this interaction is not so obvious.

\Defn\label{item-210-h}
Let $B,A$ be a pair of enriched frames with furnishings
$\JBB{B},\JDD{B}$ (on $B$) and $\JBB{A},\JDD{A}$ (on $A$),
respectively. A {\sf bed morphism} 
\OBSCURE{
\Dm B&\rTo^f& A \mD
}
\begin{center}
\includegraphics[scale=0.73,trim=0pt 0pt 0pt 0pt,clip]{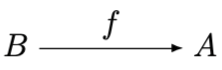}
\end{center}
is a frame morphism $f$ such that
$$(f\circ\JBB{B})(b) \leq (\JBB{A}\circ f)(b)  \hspace{5ex} 
  (f\circ\JDD{B})(b) \leq (\JDD{A}\circ f)(b)$$
for all $b\in B$. \ED
\EDefn

Notice that the requirement on a bed morphism is lax, not strict.

\section{The geometric categories}\label{sec2}

There are several geometric categories we wish to consider, all with the 
same objects, the topologically valued transition structures. These 
categories differ only in the choice of
arrows. Here we concentrate on what seems to be the most important of
these geometric categories. 

\subsection{The objects}

The basic idea is not very novel. We take a (discrete)
transition structure $(\Sblock,\rra)$ and furnish it with some
topological gadgetry. This has
been tried before 
(\cite{esakia}, \cite{gold76a}, \cite{hal55}, \cite{hil00}, 
\cite{samvac}, \cite{Wij90})
but in all cases the approach has been too simplistic. 
Usually the set $\Sblock$ of nodes is furnished with a topology, but
then to ensure that certain constructions produce sensible topologies
the transition relation is assumed to interact in some way with the
topology. A better way, which does not need such
restrictions, is to use a `fibred' setting in which the transitions
nodes are separated from the points of the space. 

\Defn
A {\sf topologically valued transition structure}
$$\Sgoth = (\Sblock,\sigma,S)$$
has three components:
\begin{itemize}
\item a transition structure $\Sblock$ with associated transition relation
$\rra$ 
\item a topological space, the {\sf value space}, with associated
  topology $\OS$ 
\item a {\sf valuation} function
\end{itemize}
\OBSCURE{
\Dm \Sblock&\rTo^{\sigma}&S \mD
}
\begin{center}
\includegraphics[scale=0.73,trim=0pt 0pt 0pt 0pt,clip]{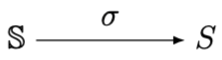}
\end{center}
which we require to be surjective. 

A {\sf plot} is just a topologically valued transition structure. \ED
\EDefn

A transition node $P\in\Sblock$ may carry lots of
quite different information. A transition
$$P\rra Q$$
in $\Sblock$ changes some, or all, of this information. Depending on
what we want to do, it may be that we need only some of the
information in a node $P$. Thus we extract this information to
give us the `value' $\sigma(P)$ of the node. We assume this value
resides in some space. In this way we have the possibility of keeping the
transition structure $\Sblock$ but changing the value space $S$ to suit our
needs. Another possibility is to use the same value space for
different transition structures. Separating these two components gives
us greater flexibility. 

We are going to talk a lot about these geometric structures, and the
official name `topologically valued transition structure' is bit
long. It is useful to have a snappy alternative name. That is the
purpose of the name `plot'.  

\Egs

(1) In many examples of plots the carrier of the transition
structure is the same as the carrier of the topological space and
the valuation is the identity function. We call these {\sf flat}
plots.

(2) Flat plots are used in various guises and under various
names to give a semantics to modal logics, see for example
\cite{esakia}, \cite{gold76a}, \cite{hil00},
\cite{Wij90}.
Several of these make use of the following construction. It is an
algebraic version of the canonical model construction found in
modal logic. Given any Boolean algebra $B$ with a unary operator
$\Box:B \rra B $ we can build a transition structure $spec (B)$ as
follows. Let the nodes be the prime filters of $B$. Define the
transition relation by:
$$P\rra Q \Longleftrightarrow \{ \Box b \mid b\in B\} \subseteq Q$$
for all prime filters $P,Q$ of $B$. The carrier set $spec (B)$ can
be endowed with the Stone topology. This has a sub-basis of opens
$$O_b = \{ P \in spec(B) \mid b \in P\}$$
where $b$ ranges across $B$.

(3) Consider a transition structure $\Sblock$ with two nodes $P$
and $Q$ and a single transition $P \rra Q$. Consider the
Sierpinski topological space $S$ with two points $P$, $Q$ where
$\{Q\}$ is open and $\{P\}$ is not open. Let $\sigma$ be the
function from $\Sblock$ to $S$ which takes $P$ to $P$ and $Q$ to
$Q$. This is an example of a flat plot in which the transition
relation coincides with the specialization order.

(4) Let the transition relation consist of the natural numbers
$\mathbb{N}$ together with the transition relation
$$P \longrightarrow Q \ \Longleftrightarrow \ P \leq Q$$
where $P,Q\in\mathbb{N}$ and $\leq$ is the standard order on
$\mathbb{N}$. Let the topology on $\mathbb{R}$ be the Euclidean
one. Let the valuation map be inclusion of the natural numbers into
the reals. In a similar way we can produce a flat plot by
considering the real numbers as a transition structure.

(5) Equilogical Spaces. According to \cite{bbs} an equilogical
space consists of a $T_0$ topological space together with a
partial equivalence relation (symmetric, transitive but not
necessarily reflexive) relation on the carrier of the space.
Generalizations of this have been considered in which the space is
not required to be $T_0$ and in which the relation may live on a
carrier other than the carrier of the space. There are also fibered
versions of equilogical spaces, in which the relation and the topology
live on different sets but are joined by a valuation from the space 
to the transition structure.

(6) Transition structures and modal logics are commonly used for the
analysis of so-called hybrid systems. A central notion is that of 
``hybrid automaton'', which consists of
a state-space for modelling the evolution of certain dynamical
systems together with a separate collection of control nodes. Many
of the standard examples of hybrid automata are in fact plots,
although some of the most general definitions 
allow the rooting map to be multi-valued.
\ED
\EEgs

\subsection{Various arrows}

Each plot is a compost of two different kinds of structures, a transition
structure and a topological space. Thus a comparison arrow for a pair of
plots should consist of a pair of
comparison arrows, one for each component. Also, there should be some
interaction with the valuation functions. We settle on the
following notion. 

\Defn\label{item-220-c}
A {\sf plot map} 
\OBSCURE{
\Dm \Sgoth&\rTo^{(\phi,\Phi)}&\Tgoth \mD
}
\begin{center}
\includegraphics[scale=0.73,trim=0pt 0pt 0pt 0pt,clip]{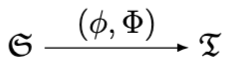}
\end{center}
from a plot
$\Sgoth=(\Sblock,\sigma,S)$ to a plot $\Tgoth=(\Tblock,\tau,T)$
is a transition morphism $\Phi$ and a continuous map $\phi$
\OBSCURE{
\Dm \Sblock&\rTo^{\Phi}&\Tblock&\hspace{3ex}& S&\rTo^{\phi}&T \mD
}
\begin{center}
\includegraphics[scale=0.73,trim=0pt 0pt 0pt 0pt,clip]{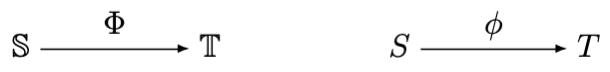}
\end{center}
such that
the square of functions on the left commutes
\OBSCURE{
\Dm[size=2em]
\Sblock       & \rTo^{\Phi}  & \Tblock    \\
\dTo^{\sigma} &             & \dTo_{\tau} &\hspace{3ex}&
(\phi\circ\sigma)(P) = (\tau\circ\Phi)(P) \\
S             & \rTo_{\phi} & T 
\mD
}
\begin{center}
\includegraphics[scale=0.73,trim=0pt 0pt 0pt 0pt,clip]{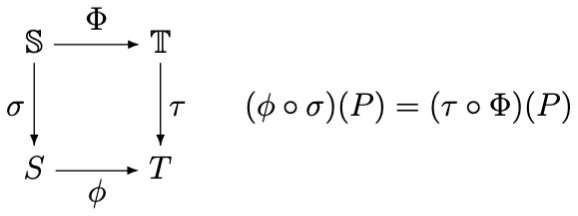}
\end{center}
which is to say that the equality on the right holds for each 
for each node $P\in\Sblock$. \ED
\EDefn

We look at quite a few plot maps, and it is useful to have a fixed
notation. Therefore Definition \ref{item-220-c} sets down the notation 
we use, as well as defining the notion.

Plot maps are composed in the obvious way by taking the function
composite of the corresponding components.  
In this way we obtain a category.

\Defn
Let $\Plot$ be the category of plots and plot maps. \ED
\EDefn

Although all of the geometric action takes place within $\Plot$, these
arrows are too weak. They are a kind of analogue of transition
morphisms. There is also an analogue for simulations, but again those
arrows are too strong. However, with our topological set-up we can
delimit several different intermediate classes of arrows, and so 
form various subcategories $\Plot^{(+)}$ with the same
objects but with a restricted family or arrows. In this paper we study
just one of these. We first set up this in some
detail, and then we indicate how other categories can be
produced by a similar method. 

\Defn
A plot map
\OBSCURE{
\Dm \Sgoth&\rTo^{(\phi,\Phi)}&\Tgoth \mD
}
\begin{center}
\includegraphics[scale=0.73,trim=0pt 0pt 0pt 0pt,clip]{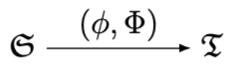}
\end{center}
is {\sf lentile} if for each node $P\in\Sblock$ and 
transition $\Phi(P)\rra R$ in $\Tblock$ we have    
$$\tau(R) \in \big((\tau\circ\Phi)[P^\ra]\big)^\ell$$
using the lens closure. \ED
\EDefn

Lentile maps are closed under composition and this gives us a 
category. 

\Defn
Let $\PlotL$ be the category of plots and lentile plot maps. \ED
\EDefn

The definition of a lentile map can be split into two
parts. We use the saturation closure
$(\cdot)^\upa$, the specialization order $\specord$, and the
topological closure $(\cdot)^-$ of the spaces.

\Lem\label{item-220-g}
A plot map $(\phi,\Phi)$, as above,  is {\sf lentile} precisely when
$(\upa)$ and $(-)$ hold.  

($\upa$) For each $P\in\Sgoth$ and situation 
$$\Phi(P)\rra R$$
in $\Tgoth$ there is some $Q\in\Sblock$ and situation 
$$P\rra Q \hspace{6ex} (\phi\circ \sigma)(Q)\specord \tau(R)$$
in $\Sgoth$ and $\Tgoth$.

($-$) For each $P\in\Sgoth$ and situation 
$$\Phi P\rra R \hspace{6ex} \tau(R)\in V\in \OT$$
in $\Tgoth$ there is some $Q\in\Sblock$ and situation 
$$P\rra Q \hspace{6ex} (\phi\circ\sigma)(Q)\in V$$
in $\Sgoth$ and $\Tgoth$.
\ELem

\Proof
Suppose first that $(\phi,\Phi)$ is lentile.

($\upa$) Consider a transition
$\Phi P\rra R$
in $\Tblock$. The lentile property gives
$$\tau R \in \big((\phi\circ\sigma)[P^\ra]\big)^\ell
    \subseteq\big((\phi\circ\sigma)[P^\ra]\big)\tsat$$
and hence there is some $Q\in\Sblock$ with
$P\rra Q$ and $(\phi\circ\sigma)(Q)\specord\tau R$
as required. 

($-$) Consider a transition
$\Phi P\rra R$ where $\tau(R)\in V\in \OT$
in $\Tgoth$. The lentile property gives
$\tau(R) \in \big((\phi\circ\sigma)[P^\ra]\big)\tclose$
so that
$(\phi\circ\sigma)[P^\ra] \mbox{ meets }V$
to give some $Q\in\Sblock$ with
$P\rra Q$ and $(\phi\circ\sigma)(Q)\in V$
as required. 

\smallskip

The converse is proved in the same way. \EP

This characterization does {\em not} say that for each
target transition $\Phi(P)\rra R$ there is a single source transition
$P\rra Q$ to witness the lentile condition. It
requires {\em two} source transitions $P\rra Q$, one for the 
$(\upa)$-condition and one for the $(-)$-condition. 

\Egs
Consider a transition morphism $\Phi:\Sblock\rTo\Tblock$
between two transition structures, and impose a topology
directly on $\Sblock$ and $\Tblock$. This gives a pair of flat
plots where the rooting maps are the identities. The only possible
continuous map
$$\phi:S=\Sblock\rTo\Tblock=T$$
is the function $\Phi$, but in general this need not be
continuous.

(General) If $\Phi$ is simulation of (discrete) transition
structures then it is also lentile. The topologies can be 
used to discard some unwanted simulations.

(Discrete) Suppose $S$ and $T$ carry the discrete
topologies. Then $\Phi$ is a plot map. By considering the
singleton sets (which are open) we see that $\Phi$ has the 
$(-)$-condition precisely when it is a simulation of (discrete)
transition structures.

(Indiscrete) Suppose $S$ and $T$ carry the indiscrete
topologies. Then $\Phi$ is a plot map. The only open sets are the
two extremes, so that $\Phi$ is a simulation precisely when
$$P^\ra=\emptyset \Rra (\Phi P)^\ra=\emptyset$$
that is, when $\Phi$ transfers dead-ends to dead-ends. Any two
points in the target space are comparable, so every plot
simulation is lentile.

(Alexandroff) Suppose $\Sblock$ and $\Tblock$ are both
partially ordered by their transition relations. So in each we
have
$$x \rra y \Longleftrightarrow x \leq y$$ for all nodes $x,y$,
Let the associated topologies be their respective Alexandroff
topologies. The continuity of the map $\Phi$ is just monotonicity.
The plot map condition is also monotonicity. The
($\uparrow$)-condition is satisfied by any function $\Phi$. The
($-$)-condition is satisfied if and only if the following
condition is satisfied:
$$\Phi (P) \leq Q \ \Rra \ (\uparrow \! P) \cap \Phi ^{\la}
(\uparrow Q) \neq \emptyset$$ for all $P \in \Sblock$ and $Q \in
\Tblock$.

($T_1$) Suppose the target set $\Tblock$ carries a $T_1$
topology. Then an arrow with the ($\uparrow$)-condition 
must be a (discrete) simulation. However, not every (discrete) 
simulation need be a plot map.

($\mathbb{R}$) Consider the flat plot with base space the reals
$\mathbb{R}$ and transition relation given by the order on the
reals:
$$x \longrightarrow y \ \Longleftrightarrow \ x \leq y$$ for all
$x,y \in \mathbb{R}$.. Consider the function $\Phi:\mathbb{R}
\rightarrow \mathbb{R} ; x \mapsto x$. This is  lentile map from
the plot to itself. If a plot map satisfies either the
$(-)$-condition or the ($\uparrow$)-condition then it is a
discrete simulation. On the other hand, the constant map
$x\mapsto 0$ is an example of a plot map which does not satisfy
the $(-)$-condition.

(A plot map and homeomorphism need not satisfy the
$(-)$-condition.)
Let $\Sblock$ be the flat plot consisting of the Sierpinski space
$\mathbf{2}$ together with the empty transition relation. Let
$\Tblock$ be the flat plot consisting of the space $\mathbf{2}$
together with the transition relation $\longrightarrow$ with
$$(x\longrightarrow y \ \Longleftrightarrow \ x=0 \ \& \ y=1)$$
for all $x,y$. Then the map $\Phi:X \rightarrow Y$ with
$\Phi(0)=0$ and $\Phi(1)=1$ is a homeomorphism and a plot map but
does not satisfy the $(-)$-condition.

(A lentile homeomorphism.) 
Consider flat plots $(\Sblock,S)$ and $(\Tblock,T)$. Suppose
$\Phi:S \rightarrow T$ is a homeomorphism, with inverse $\Phi
^{-1}$. If $\Phi$ and $ \Phi ^{-1}$ are both plot maps then $\Phi$
and $\Phi^{-1}$ are both lentile.

(Plot map which satisfies the ($\uparrow$)-condition  but is 
not induced by a simulation.)
Let the transition structure $\Sblock$ have one node $P$ and one
transition. Let the transition structure $\Tblock$ have two nodes,
$Q$ and $R$. Let the transitions be $Q \rra Q$ and $Q \rra R$. Let
$S$ be the one-point space and call its unique element $p$. Define
plots $(\Sblock,\sigma,S)$ and $(\Tblock,\tau,S)$ with rooting
maps $\sigma$ and $\tau$. Consider
the continuous map $\phi :S \rra S$ as the identity. We define the
transition morphism component $\Phi: \Sblock \rra \Tblock$ by
letting $\Phi P =Q$. The plot morphism square for $(\phi,\Phi)$
commutes since $\Sblock$ has only one point $P$ and
$$( \phi \circ \sigma )(P) =\phi (\sigma (P)) = \phi (p) = p =
\tau (Q) = \tau ( \Phi (P)) = ( \tau \circ \Phi )(P)$$ holds. The
plot map is tight since $$(\phi \circ \sigma)(P)=\tau Q =\tau R$$
and $P \longrightarrow P$ hold. The map $\Phi$ is not a
simulation, for $\Phi P$ makes a transition to $R$ but there is no
point which $\Phi$ maps to $R$.
\ED
\EEgs

For much of what we do the base space is fixed. Thus we have some
space $S$ and we consider different plots all valued over
$S$. The comparison maps between these all have the identity
functions on $S$ as the spatial component. In other words, we work in
$\Plot_S$, the {\sf $S$-based fibre} of $\Plot$. (We do not
claim that the functor from $\Plot$ is to spaces is a fibration.)
In some ways this fibred view of $\Plot$ give a better perspective. In
particular, an analysis of $\Plot_S$ is somehow connected with the 
various enrichments $\OS$ can carry.

To conclude this subsection we say a few words about the more general
construction. The trick is to replace the lens closure $(\cdot)^\ell$
by some other suitable closure operation $(\cdot)^+$. For instance
both $(\cdot)^\upa$ and $(\cdot)^-$ are suitable, and there are
others. We then consider those plot maps such that 
$$\tau(R) \in \big((\tau\circ\Phi)[P^\ra]\big)^+$$
for each target transition $\phi P\rra R$. This gives us several
categories $\Plot^{(+)}$ whose properties are studied elsewhere. 

\subsection{Enriching the topology}\label{subsec2.3}

Each plot $\Sgoth$ has a discrete transition structure $\Sblock$ as one
component. In the usual way the transition relation on $\Sblock$
induces operations $\JBB{\Sblock},\JDD{\Sblock}$ on
$\Pcal\Sblock$. The other component of $\Sgoth$ is a space $S$ with
its topology $\OS$. These components are connected by the valuation
function $\sigma$. We might expect that $\sigma$ interacts in some
with way the discrete furnishings of $\Pcal\Sblock$, except, as yet,
there are no furnishings at the other end of $\sigma$. We shall 
rectify that. 

The valuation function $\sigma$ has an inverse image function
\OBSCURE{
\Dm \OS&\rTo^{\sigma^\la}&\Pcal\Sblock \mD
}
\begin{center}
\includegraphics[scale=0.73,trim=0pt 0pt 0pt 0pt,clip]{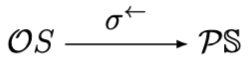}
\end{center}
which is a frame morphism. In other words, $\sigma$ is a continuous
map when we remember that $\Pcal\Sblock$ is the discrete topology on
$\Sblock$. We furnish $\OS$ so that $\sigma^\la$ becomes an
enriched frame morphism. To do that we merely lift the furnishings on
$\Pcal\Sblock$ across $\sigma^\la$.

\Defn\label{item-220-i}
Let $\Sgoth=(\Sblock,\sigma,S)$ be a plot. The {\sf lift} of the
operations $\JBB{\Sblock},\JDD{\Sblock}$ on $\Pcal\Sblock$ are the
operations $\JBB{\sigma},\JDD{\sigma}$ on $\OS$ given by
$$V\subseteq\JBB{\sigma}(U) \LRa 
\sigma^\la(V)\subseteq(\JBB{\Sblock}\circ\sigma^\la)(U)  \hspace{5ex}
  V\subseteq\JDD{\sigma}(U) \LRa 
\sigma^\la(V)\subseteq(\JDD{\Sblock}\circ\sigma^\la)(U)$$
for $U,V\in\OS$. \ED
\EDefn

This construction drives everything we do. 
Let's set down its basic properties.

\Theom\label{item-220-j}
For each plot $\Sgoth=(\Sblock,\sigma,S)$ the lifted furnishings convert
$\OS$ into a bed in such a way that $\sigma^\la$ is a bed morphism. 
\ETheom

\Proof
In more detail we are claiming that 
$$\JBB{\sigma}(S)=S \hspace{5ex} \JDD{\sigma}(\emptyset)= \emptyset$$
and
$$\JBB{\sigma}(U\cap V) = \JBB{\sigma}(U) \cap \JBB{\sigma}(V) 
\hspace{5ex} V\subseteq U \Rra \JDD{\sigma}(V)\subseteq\JDD{\sigma}(U)$$
together with the mixed condition
$$\JBB{\sigma}(U) \cap \JDD{\sigma}(V) \subseteq \JDD{\sigma}(U\cap V)$$
for $U,V\in\OS$. Only the mixed condition is not immediate, 
but this is still easy. Surjectivity of $\sigma$ is used in showing $\JDD{\sigma}(\emptyset) \subseteq \emptyset$.

For the bed morphism we require
$$\sigma^\la\circ\JBB{\sigma} \leq \JBB{\Sblock}\circ\sigma^\la
\hspace{5ex}
  \sigma^\la\circ\JDD{\sigma} \leq \JDD{\Sblock}\circ\sigma^\la$$
but these are essentially the definitions of $\JBB{\sigma}$ and
$\JDD{\sigma}$. \ED

\section{The algebraic categories}\label{sec3}

We refine the notion of a bed by attaching a space to each
and insisting that the comparison morphisms between 
these new gadgets also has some topological content. 

\subsection{The objects}
In the same way that we moved from a transition structure to a plot,
we move from a bed to a more refined object. 

\Defn
A {\sf garden} 
$$\Agoth = (A,\alpha,S)$$
over a space $S$ has three components
\begin{itemize} 

\item a bed $A$ 
\item a topological space $S$ 
\item a {\sf covering} morphism
\OBSCURE{
\Dm A&\rTo^\alpha&\OS \mD
}
\begin{center}
\includegraphics[scale=0.75,trim=0pt 0pt 0pt 0pt,clip]{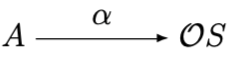}
\end{center}
\end{itemize}
which is a surjective frame morphism. \ED
\EDefn

Notice that the covering morphism $\alpha$ is merely a frame morphism,
and does not interact with the furnishings. It can't, for as yet the
topology carries no furnishings. 

Notice also that the covering morphism must be surjective. This is
important and we use the condition at one or two crucial
places. We point out all these uses.

\Eg\label{item-230-b} 
Each bed $A$ gives a garden. 
We take the point space $S$ of $A$ as the base space 
and the canonical morphism $A\rTo\OS$ as the
covering morphism. \ED
\EEg

The next example has the status of a minor result. It is 
a rephrasing of Theorem \ref{item-220-j}.

\Lem\label{item-230-c}
For each plot 
$$\Sgoth=(\Sblock,\sigma,S)$$ 
the induced operations on the topology $\OS$ produce a garden
$$(\OS,\id,S)$$
over the same space with a trivial covering morphism.
\ELem

For both these examples the base space $S$ has a rather special
connection with the bed $A$. In fact, there is a temptation
to consider only those gardens of Example \ref{item-230-b}. However, by
uncoupling the space from the frame we obtain a bit more flexibility.

Notice that a plot $\Sgoth=(\Sblock,\sigma,S)$ does {\em not} give a
garden of the form
$(\Pcal\Sblock,\alpha,S)$
for, in general, there is no sensible frame morphism
$\Pcal\Sblock\rTo\OS$
which can be used as the covering morphism. In particular, the
interior operation $(\cdot)^\circ$ is {\em not} a frame morphism. 

\subsection{The arrows}

As with a plot, each garden is a compost of two different kinds of
structures. Thus we need a pair of comparison arrows, one for each 
component, and there should be some interaction with the covering 
morphism. The actual notion we use is not quite so obvious. 
In the following definition we set down some standard notation.

\Defn\label{item-230-d}
A {\sf garden morphism}
\OBSCURE{
\Dm \Bgoth&\rTo^{(\phi,f)}&\Agoth \mD
}
\begin{center}
\includegraphics[scale=0.73,trim=0pt 0pt 0pt 0pt,clip]{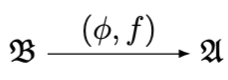}
\end{center}
from a garden 
$\Bgoth=(B,\beta,T)$ to a garden $\Agoth=(A,\alpha,S)$
is a bed morphism $f$ and a
continuous map $\phi$
\OBSCURE{
\Dm B&\rTo^{f}&A&\hspace{3ex}& T&\lTo^{\phi}&S \mD
}
\begin{center}
\includegraphics[scale=0.73,trim=0pt 0pt 0pt 0pt,clip]{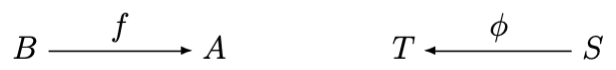}
\end{center}
such that
the square of frame morphisms on the left commutes
\OBSCURE{
\Dm[size=2em]
       B     & \rTo^{f}        & A \\
\dTo^{\beta} &                 & \dTo_{\alpha} &\hspace{3ex}&
(\phi^\la\circ\beta)(b) = (\alpha\circ f)(b)     \\
\OT            & \rTo_{\phi^\la} & \OS 
\mD
}
\begin{center}
\includegraphics[scale=0.73,trim=0pt 0pt 0pt 0pt,clip]{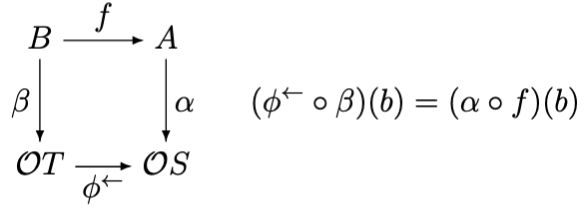}
\end{center}
that is the equality on the right holds for each 
element $b\in B$. \ED
\EDefn

There are a couple of things to note. The continuous map $\phi$
seems to go in the wrong direction. However, we use the inverse
image function $\phi^\la$ as a frame morphism and this does go in the
right direction.  The interaction condition is strict and not a lax
comparison. 

\Defn
Let $\Grdn$ be the category of gardens and garden morphisms. \ED
\EDefn

Just as with $\Plot$ we may consider the fibre category
$\Grdn_S$
of gardens over a fixed base space $S$. In such a category the
spatial component of each morphism is the identity function on $S$. In
many ways this gives a better perspective of the various
constructions. In particular, each of the functors we produce can be
restricted to operate over $S$. (As with $\Plot$, we are not claiming
that the functor from $\Grdn$ to spaces is a fibration.)

\subsection{The functor}

Lemma \ref{item-230-c} shows that each plot $\Sgoth$ gives us a garden
over the same space. There is a companion construction on arrows.

Consider a pair of plots
$$\Sgoth=(\Sblock,\sigma,S) \hspace{5ex}\Tgoth=(\Tblock,\tau,T)$$
with induced operations
$$\JBB{\sigma}\hspace{3ex}\JDD{\sigma}\hspace{7ex}
  \JBB{\tau}\hspace{3ex}\JDD{\tau}$$
on the respective topologies $\OS$ and $\OT$. 
Consider also a map 
\OBSCURE{
\Dm \Sgoth&\rTo^{(\phi,\Phi)}&\Tgoth \mD
}
\begin{center}
\includegraphics[scale=0.73,trim=0pt 0pt 0pt 0pt,clip]{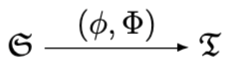}
\end{center}
between plots. The function $\phi$
is a continuous map between the spaces, and so gives a
frame morphism $\phi^\la$
between the topologies. Thus we have the components
\OBSCURE{
\Dm \OT&\rTo^{\phi^\la}&\OS &\hspace{3ex}&T&\lTo^{\phi}&S \mD
}
\begin{center}
\includegraphics[scale=0.73,trim=0pt 0pt 0pt 0pt,clip]{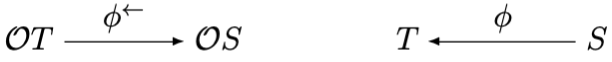}
\end{center}
required for a garden morphism. 
For $\phi^\la$ to be a bed morphism we must start with
the right kind of plot map.

\Lem
If a plot map $(\phi,\Phi)$ is lentile then $(\phi,\phi^\la)$ is a
$\Grdn$-arrow.  
\ELem

\Proof
We show that the two comparisons
$$\phi^\la\circ\JBB{\tau} \leq \JBB{\sigma}\circ\phi^\la
\hspace{10ex}
  \phi^\la\circ\JDD{\tau} \leq \JDD{\sigma}\circ\phi^\la$$
hold. The initial parts of the two verifications are similar so we 
do these in parallel. In other words, in the next paragraph, the 
left-hand displays and the right-hand displays should be read independently. 

Consider any $W\in OT$ and let
$$V=(\phi^\la\circ\JBB{\tau})(W) \hspace{3ex} U=\phi^\la(W)
\hspace{7ex}
  V=(\phi^\la\circ\JDD{\tau}(W) \hspace{3ex} U=\phi^\la(W)$$
to obtain $U,V\in OS$. We require
$$V\subseteq \JBB{\sigma}(U) \hspace{10ex} V\subseteq \JDD{\sigma}(U)$$
so, by Definition \ref{item-220-i}, an  inclusion
$$\sigma^\la(V) \subseteq (\JBB{S}\circ\sigma^\la)(U)
\hspace{7ex}
  \sigma^\la(V) \subseteq (\JDD{S}\circ\sigma^\la)(U)$$
will suffice. We look at these separately.

For the left-hand inclusion, using the definition of $V$, the plot map
property, the construction of $\JBB{\tau}$, the discrete morphism
property, and the plot map property again, and finally the definition
of $U$ we have  
$$\begin{array}{lclcl}
             \sigma^\la(V)
&=        & (\sigma^\la\circ\phi^\la\circ\JBB{\tau})(W) \\
&=        & (\Phi^\la\circ\tau^\la\circ\JBB{\tau})(W)   \\
&\subseteq& (\Phi^\la\circ\JBB{\Tblock}\circ\tau^\la)(W)   \\
&\subseteq& (\JBB{\Sblock}\circ\Phi^\la\circ\tau^\la)(W)   \\
&=&         (\JBB{\Sblock}\circ\sigma^\la\circ\phi^\la)(W)   
&=&         (\JBB{\Sblock}\circ\sigma^\la)(U)   
\end{array}$$
as required. For this part none of the lentile properties
are used.

For the right-hand inclusion, using the definition of $V$, 
the plot map property,  the construction of $\JBB{\tau}$, the
$(-)$-property of a lentile map, the plot map property again, and
finally the definition of $U$ we have
$$\begin{array}{lclcl}
             \sigma^\la(V)
&=        & (\sigma^\la\circ\phi^\la\circ\JDD{\tau})(W) \\
&=        & (\Phi^\la\circ\tau^\la\circ\JDD{\tau})(W)   \\
&\subseteq& (\Phi^\la\circ\JDD{T}\circ\tau^\la)(W)   \\
&=&         (\JDD{S}\circ\sigma^\la\circ\phi^\la)(W)   
&=&         (\JDD{S}\circ\sigma^\la)(U)   
\end{array}$$
as required. For this part some but not all of the lentile
properties are used. \EP

This more or less completes the proof of the following. 

\Theom
There is a contravariant functor
\OBSCURE{
\Dm \PlotL&\rTo^{\Gcal}&\Grdn \mD
}
\begin{center}
\includegraphics[scale=0.73,trim=0pt 0pt 0pt 0pt,clip]{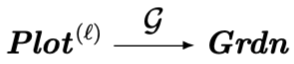}
\end{center}
which sends each plot into a furnished version of the topology of the
value space.
\ETheom

So far this is not very exciting, and if this was all we had to say
then it wouldn't be worth saying. The main interest is to find
a contravariant adjoint to $\Gcal$. The remainder of
the paper is devoted to that functor and its properties.

\section{Harvesting a garden}\label{sec4}

The functor $\Gcal$ converts a plot into a garden. The
gardens produced in this way are rather special since they are
furnished versions of the topology of the base space. 
We now begin the process of going in the other direction. We first produce
a contravariant functor
\OBSCURE{
\Dm[width=4em] \Grdn&\rTo^{\Fcal} &\PlotL \mD
}
\begin{center}
\includegraphics[scale=0.73,trim=0pt 0pt 0pt 0pt,clip]{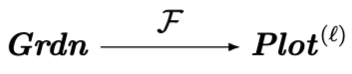}
\end{center}
and then we show that the pair $\Fcal,\Gcal$ form a
contravariant adjunction. 

\Defn\label{item-240-a}
(a) For a frame $A$ let $\Fil(A)$ be the set of filters on $A$.

(b) Let $\Agoth=(A,\alpha,S)$ be a garden over a space $S$. For $p\in
S$ we use 
$$x\in\nabla(p)  \LRa p\in\alpha(x)            \hspace{2.5ex}
  x\in\PBB{A}(p) \LRa \JBB{A}(x)\in\nabla(p)   \hspace{2.5ex}
  x\in\PDD{A}(p) \LRa \JDD{A}(x)\notin\nabla(p)$$
(for $x\in A$) to produce three assignments
$$\nabla,\PBB{A},\PDD{A}:S\rTo\Pcal A$$
from points of $S$ to subsets of $A$. \ED
\EDefn

Note that $\nabla(p),\PBB{A}(p)\in\Fil(A)$, and 
$\PDD{A}(p)$ is a lower section of $A$.

\subsection{The object assignment}

To convert a garden $\Agoth=(A,\alpha,S)$ into a cultivated
plot $\Sgoth=(\Sblock,\sigma,S)$ over the same space $S$
we set up a suitable transition structure $\Sblock$.
For this, each transition node is a triple where
the first component is a point $p\in S$. The valuation function
\OBSCURE{
\Dm \Sblock&\rTo^{\proj}&S \mD
}
\begin{center}
\includegraphics[scale=0.73,trim=0pt 0pt 0pt 0pt,clip]{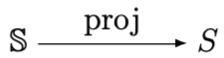}
\end{center}
is then the projection which selects the first component. 

\Defn\label{item-240-b}
Let $\Agoth=(A,\alpha,S)$ be a garden over a space $S$. A 
{\sf flower} of $\Agoth$ is a triple
$$(p,a,F)$$
where the {\sf root} $p\in S$, the {\sf stalk} $a\in A$, and  
the {\sf bloom} $F\in\Fil(A)$
must satisfy
$$a\in\PDD{A}(p)   \hspace{5ex} \PBB{A}(p)\subseteq F$$
using the associated operators $\PDD{A},\PBB{A}$ of Definition
\ref{item-240-a}. 

Let $\Fblock$  be the set of flowers of $\Agoth$. \ED
\EDefn

The conditions on a flower can be rephrased. We have
$$a\in\PDD{A}(p) \LRa \JDD{A}a\notin \nabla(p) 
                 \LRa p\notin(\alpha\circ\JDD{A})(a)$$
and
$$\PBB{A}(p)\subseteq F \LRa 
(\forall x)\big[p\in (\alpha\circ\JBB{A})(x) 
\Rra \JBB{A}(x)\in\nabla(p) \Rra x\in F\big]$$
both of which we use at various times. 

This attaches to $\Agoth$ a set of nodes $\Fblock$. To 
complete the object construction we do two things to 
this set. We give it a transition relation,
and then we extract a certain subset. 

\Defn\label{item-240-c}
The set $\Fblock$ of flowers of a garden carries a transition relation
given by  
$$(p,a,F)\rra (q,b,G) \LRa 
a\notin\nabla(q)\mbox{ and }F\subseteq\nabla(q) $$
for flowers $(p,a,F),(q,b,G)\in\Fblock$. \ED
\EDefn

This converts the set $\Fblock$ into a transition structure. 
However, in general this set is too big and some of the
flowers are not in good
condition. We have to weed out the deadwood.

\Defn
A set $\Hblock\subseteq\Fblock$ of flowers of a garden is 
{\sf healthy} if both 
\begin{itemize}
    \item[$(\upa)$] For each $(p,a,F)\in\Hblock$ and each $x\notin F$
      there is a transition  
      $$(p,a,F)\rra(q,b,G)$$ 
      with  $(q,b,G)\in\Hblock$ and $x\notin\nabla(q)$

    \item[$(-)$] For each $(p,a,F)\in\Hblock$ and each $x\not\leq a$
     there is a transition  
     $$(p,a,F)\rra(q,b,G)$$
     with  $(q,b,G)\in\Hblock$ and $x\in\nabla(q)$
\end{itemize}
hold. A flower is healthy if it belongs to some healthy set. \ED
\EDefn

The union of a family of healthy sets of flowers
(from the same garden) is itself healthy. Thus each garden
$\Agoth$ has a largest healthy subset of $\Fblock$. This is 
the important set. 

\Defn\label{item-240-e}
For each garden $\Agoth = (A,\alpha,S)$ let $\Sblock$ be the
largest healthy subset of $\Fblock$ with the transition relation
induced from $\Fblock$. Let 
$$\Fcal\Agoth = (\Sblock,\proj,S)$$
where the valuation function is the projection to the base space $S$
of points. \ED 
\EDefn

\subsection{Passing across a morphism}\label{subsec4.2}

To construct the arrows assignment we fix some notation. Let
\OBSCURE{
\Dm \Bgoth&\rTo^{(\phi,f)}&\Agoth \mD
}
\begin{center}
\includegraphics[scale=0.73,trim=0pt 0pt 0pt 0pt,clip]{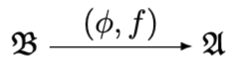}
\end{center}
be a garden morphism between gardens
$$\Bgoth=(B,\beta,T) \hspace{5ex} \Agoth=(A,\alpha,S)$$
with the usual furnishings. 
Notice that we do not assume that the base 
spaces are the same.  
We continue with the notation set up in
Definition \ref{item-230-d}. Recall also Definition \ref{item-210-h}.

\Lem\label{item-240-f}
For the garden morphisms $(\phi,f)$, as above, we have
$$(\nabla\circ\phi)(p)=(f^\la\circ\nabla)(p) \hspace{5ex}
(\PBB{B}\circ\phi)(p) \subseteq (f^\la\circ\PBB{A})(p)$$
for each $p\in S$. 
\ELem

\Proof
For each point $p\in S$ and element $y\in B$ we have
$$y\in(\nabla\circ\phi)(p)
\LRa \phi(p)\in\beta(y)
\LRa p\in(\phi^\la\circ\beta)(y)=(\alpha\circ f)(y)
\LRa f(y)\in\nabla(p)$$
Similarly we have
$$y\in\big(\PBB{B}\circ\phi\big)(p) 
\LRa \JBB{B}(y) \in \big(\nabla\circ\phi\big)(p)
=   \big(f^\la\circ\nabla\big)(p)
\LRa \big(f\circ\JBB{B}\big)(y) \in \nabla(p)$$
using the first part, and then the morphism property gives
$$y\in\big(\PBB{B}\circ\phi\big)(p) 
\Rra \big(\JBB{A}\circ f\big)(y) \in \nabla(p)
\Rra f(y) \in \PBB{A}(p)$$
which leads to the right-hand result. \EP

Our aim is to produce a plot map
\OBSCURE{
\Dm \Fcal\Agoth&\rTo^{\hspace{6ex}} \Fcal\Bgoth \mD
}
\begin{center}
\includegraphics[scale=0.73,trim=0pt 0pt 0pt 0pt,clip]{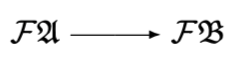}
\end{center}
so we require a transition morphism between
the two sets of healthy flowers. To get at this we first set up
a transition morphism between the full sets of (not necessarily
healthy) flowers $\Gblock$ for $\Bgoth$ and $\Fblock$ for $\Agoth$.

Recall that, since $f$ is a frame morphism, it has a right
adjoint $f_\ast$ (which transfers elements of $A$ to elements of $B$),
and an inverse function $f^\la$ which certainly transfers filters on
$A$ to filters on $B$. We use these to construct the transition
morphism.

\Lem
For the situation above we have
$$\big(\phi(p),f_\ast(a),f^\la(F)\big) \in \Gblock$$
for each flower $(p,a,F)\in\Fblock$.
\ELem

\Proof
We are given a flower $(p,a,F)\in\Fblock$ as in Definition
\ref{item-240-c}. Trivially we have
$$\phi(p)\in T \hspace{5ex} f_\ast(a)\in B \hspace{5ex}
f^\la(F)\in\Fil(B)$$
and we must check that
$$f_\ast(a)\in \big(\PDD{B}\circ\phi\big)(p) \hspace{7ex} 
\big(\PBB{B}\circ\phi\big)(p)\subseteq f^\la(F)$$
hold.

For the left hand requirement we observe that
$$\phi^\la\circ\beta\circ\JDD{B}\circ f_\ast
=    \alpha\circ f\circ\JDD{B}\circ f_\ast
\leq \alpha\circ\JDD{A}\circ f\circ f_\ast 
\leq \alpha\circ\JDD{A}$$
holds. The first step follows since
$(\phi,f)$ is a garden morphism, the second since $f$ is a bed 
morphism, and the third since the composite $f\circ f_\ast$ is
deflationary. 
But now
$$\begin{array}{l@{\;\Rra\;}lc@{\Rra\;}l}
  f_\ast(a)\notin \big(\PDD{B}\circ\phi\big)(p)      &
  \big(\JDD{B}\circ f_\ast\big)(a) \in \big(\nabla\circ\phi\big)(p) \\[0.5ex]
& \phi(p)\in (\beta\circ\JDD{B}\circ f_\ast)(a)         \\[0.5ex]
& p\in (\phi^\la\circ\beta\circ\JDD{B}\circ f_\ast)(a) \\[0.5ex]
& p \in (\alpha\circ\JDD{A})(a)                        \\[0.5ex]
& \JDD{A}(a) \in \nabla(p)                           &
& a\notin\PDD{A}(p) 
\end{array}$$
which, by taking the contrapositive, leads to the required result 

For the right hand inclusion by a use of Lemma \ref{item-240-f} we have
$$y\in\big(\PBB{B}\circ\phi\big)(p)  \Rra f(y) \in \JBB{A}(p) \subseteq F$$
to give the required result. \EP

This result shows that the garden morphism $(\phi,f)$ induces a
commuting square 
\OBSCURE{
\Dm[height=1.5em,width=3em]
\Fblock & \rTo^\Psi             & \Gblock       \\
\dTo    &                       & \dTo          &\hspace{2ex}&
(p,a,F) &\rMapsto^{\hspace{4ex}}& (\phi(p),f_\ast(a),f^\la(F)) \\
  S     & \rTo_{\phi}           & T
\mD
}
\begin{center}
\includegraphics[scale=0.73,trim=0pt 0pt 0pt 0pt,clip]{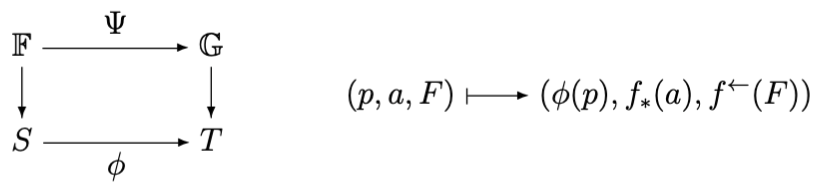}
\end{center}
of functions, as on the left. The vertical arrows are the
projections from flowers to roots. The top
function $\Psi$ is described on the right.
We show this is a transition morphism. 

\Lem\label{item-240-h}
For the situation above we have
$$(p,a,F)\rra (q,b,G)  \Rra 
\big(\phi(p),f_\ast(a),f^\la(F)\big)\rra\big(\phi(q),f_\ast(b),f^\la(G)\big)$$
for each transition from $\Fblock$.
\ELem

\Proof
By Definition \ref{item-240-c} it suffices to show that
$$f_\ast(a)\in \big(\nabla\circ\phi\big)(q) \Rra a\in\nabla(p) 
\hspace{5ex}
F\subseteq \nabla(q) \Rra f^\la(F)\subseteq\big(\nabla\circ \phi\big)(q)$$
hold. 

For the left hand implication, assuming
$$f_\ast(a) \in \big(\nabla\circ\phi\big)(q) 
= \big(f^\la\circ\nabla\big)(q)$$
we have
$$(f\circ f_\ast)(a) \in\nabla(q)$$
and hence $a\in\nabla(q)$ since $f\circ f_\ast$ is deflationary and
$\nabla(q)$ is a filter.

For the right hand implication assuming $F\subseteq\nabla(q)$ we have
$$f^\la(F) \subseteq \big(f^\la\circ\nabla\big)(q) 
= \big(\nabla\circ\phi\big)(q)$$
which is what we want. \EP

Of course, we are not interested in all flowers,
only the healthy ones. Thus we need to restrict the actions of the 
assignment $\Psi$. 

\Lem\label{item-240-i}
If $\Hblock\subseteq\Fblock$ is healthy (in $\Agoth$) then
$\Psi[\Hblock]\subseteq\Gblock$ is healthy (in $\Bgoth$).  
\ELem

\Proof
Consider a typical member of $\Psi[\Gblock]$. This has the form
$$\big(\phi(p),f_\ast(a), f^\la(F)\big)$$
for some $(p,a,F)\in\Hblock$. We must produce two transitions
$$\big(\phi(p),f_\ast(a), f^\la(F)\big) \rra 
  \big(\phi(q),f_\ast(b),f^\la(G)\big)$$
in $\Psi[\Hblock]$ subject to certain conditions. We use the
healthiness of $\Hblock$ to produce transitions
$$(p,a, F) \rra (q,b,G)$$
in $\Hblock$ and then apply Lemmas \ref{item-240-h} and \ref{item-240-f}.

$(\upa)$ Consider any element $y$ of $B$ with $y\notin f^\la(F)$. Thus
$f(y)\notin F$ and so the $(\upa$)-healthiness of $\Hblock$ gives a transition,
as above, with $f(y)\notin\nabla(q)$. But now
$$y \notin \big(f^\la\circ\nabla\big)(q) 
= \big(\nabla\circ\phi\big)(q)$$
to give the required result. 

$(-)$ Consider any element $y$ of $B$ with $y\nleq f_\ast(a)$. Thus
$f(y)\nleq a$ and so the $(-)$-healthiness of $\Hblock$ gives a transition,
as above, with $f(y)\in\nabla(q)$. But now
$$y \in \big(f^\la\circ\nabla\big)(q) = \big(\nabla\circ\phi\big)(q)$$
to give the required result. \EP

This completes the construction of the functor $\Fcal$. The
object assignment is given by Definition \ref{item-240-e}. The
arrow assignment is the restriction of the function
$\Psi$, as justified by Lemma \ref{item-240-i}. 
The checks required to show that $\Fcal$ is functorial are
straightforward.

\Defn
For each garden $\Agoth = (A,\alpha,S)$ we set
$$\Fcal\Agoth = (\Sblock,\proj,S)$$
to obtain a plot. Here $\Sblock$ is the
maximum healthy set of flowers and $\proj$ is the projection.

For each garden morphism $(\phi,f)$ we set
$$\Fcal(\phi,f) = (\phi,\Psi):\Fcal\Agoth\rTo\Fcal\Bgoth$$
to obtain a plot map. \ED
\EDefn

\subsection{Functorial properties}\label{subsec4.3}

In the previous subsection we deliberately didn't
mention the target category of the constructed functor
$\Fcal$. Strictly speaking, all that we have done so far is to produce
a  functor to the category $\Plot$ of plots and all plot
maps, but we can do better than that.

The proof of the next result uses the surjectivity of the covering
morphism. 

\Lem
For each garden morphism
\OBSCURE{
$\Dm[width=2em] \Bgoth&\rTo^{(\phi,f)}&\Agoth \mD$
}
\includegraphics[scale=0.73,trim=1ex 0.7ex 0.5ex 1ex,clip]{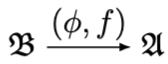}
the induced plot map
$\Fcal\Agoth\rTo\Fcal\Bgoth$
is a lentile map.
\ELem

\Proof
As usual the proof comes in two parts using
the characterization of Lemma \ref{item-220-g}. For both part we start from a
flower $(p,a,F)\in\Sblock$ and a transition
$$\big(\phi(p),f_\ast(a), f^\la(F)\big) \rra (r,c,H)$$
in $\Tblock$ (with an extra condition in the $(-)$-case), and we
must produce a transition
$$(p,a, F) \rra (q,b,G)$$
in $\Sblock$ satisfying a certain condition. 

\smallskip

$(\upa$) Given the $\Tblock$-transition, as above, we require an
$\Sblock$-transition, as above, with  $\phi(q)\specord r$ (in $T$). 
Consider the closed set $r\tclose$ of $T$, Since the covering map
$\beta$ is surjective we have
$$r\tclose = \beta(y)'$$
for some $y\in B$. In particular, $y\notin\nabla(r)$. 
We use this with
the $(\upa)$-healthiness of $\Sblock$ to produce $(q,b,G)$.

From the given $\Tblock$-transition we have 
$f^\la(F)\subseteq\nabla(r)$ and hence $f(y)\notin F$ (otherwise
$y\in\nabla(r)$, which is not so). The healthiness gives a
$\Sblock$-transition, as above, with $f(y)\notin\nabla(q)$. 
Thus, by Lemma \ref{item-240-f}, we have
$$y\notin\big(f^\la\circ\nabla\big)(q)=\big(\nabla\circ\phi\big)(q)$$
and hence $\phi(q)\in \beta(y)'=r\tclose$, to give $\phi(q)\specord r$,
as required.

$(-)$ We are given a $\Tblock$-transition, as above, together with
some $V\in\OT$ where $r\in V$. We must produce an $\Sblock$-transition,
as above, with $\phi(q)\in V$. 
We use the $(-)$-healthiness of $\Sblock$ to produce the
$\Sblock$-transition. 

From the given $\Tblock$-transition we have 
$f_\ast(a)\notin\nabla(r)$. 
The cover $\beta$ is surjective, so $V=\beta(y)$ for some $y\in B$. Then
$$f(y)\leq a 
\Rra y\leq f_\ast(a)
\Rra y\notin\nabla(r) 
\Rra r\notin\beta(y) = V$$ 
and hence $fy\nleq a$. 
The $(-)$-healthiness of $\Sblock$ gives a $\Sblock$-transition, as
above, with $f(y)\in\nabla(q)$. Thus, by Lemma \ref{item-240-f}, we have
$$y\in\big(f^\la\circ\nabla\big)(q)=\big(\nabla\circ\phi\big)(q)$$
and hence $\phi(q)\in \beta(y)=V$, to give the required result. \EP

With this we have the result we are looking for.

\Theom
There is a contravariant functor
\OBSCURE{
\Dm[width=4em] \Grdn&\rTo^\Fcal &\PlotL \mD
}
\begin{center}
\includegraphics[scale=0.73,trim=0pt 0pt 0pt 0pt,clip]{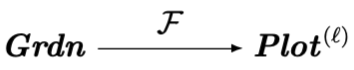}
\end{center}
formed by harvesting the garden.
\ETheom

\section{Adjunction properties}\label{sec5}

We show that the pair of functors
\OBSCURE{
\Dm[width=4em]
\Grdn &\pile{\rTo^{\Fcal} \\ \\ \lTo_{\Gcal}} &\PlotL
\mD
}
\begin{center}
\includegraphics[scale=0.73,trim=0pt 0pt 0pt 0pt,clip]{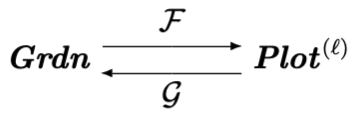}
\end{center}
form a contravariant adjunction.
We analyse the units
\OBSCURE{
\Dm
\Agoth&\rTo^{\AunitA}&(\Gcal\circ\Fcal)\Agoth&&
\Sgoth&\rTo^{\GunitS}&(\Fcal\circ\Gcal)\Sgoth
\mD
}
\begin{center}
\includegraphics[scale=0.73,trim=0pt 0pt 0pt 0pt,clip]{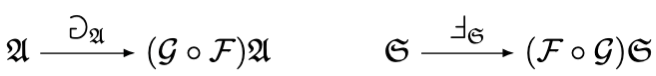}
\end{center}
in some detail. We show that each is a natural transformation, and the
transpose  
$$\Fcal(\AunitA) \hspace{5ex} \Gcal(\GunitS)$$
of each is an isomorphism, and so the adjunction is idempotent.

\subsection{The algebraic unit}\label{subsec5.1}

For the endo-functor
$$\Mcal=\Gcal\circ\Fcal$$
on $\Grdn$, we attach to each garden $\Agoth$ a garden morphism.
\OBSCURE{
\Dm \Agoth&\rTo^{\Aunit_\Agoth}&\Mcal\Agoth \mD
}
\begin{center}
\includegraphics[scale=0.73,trim=0pt 0pt 0pt 0pt,clip]{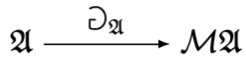}
\end{center}
To do this we fix some
notation, as in Table \ref{table-alg}, and review the
relevant constructions.

We start with a garden $\Agoth$ as in the left-hand column of the
table. The frame $A$ carries furnishings $\JBB{A},\JDD{A}$.
The functor $\Fcal$ converts this garden into a plot, as in the second
column. The transition nodes, the healthy flowers, are certain triples
$(p,a,F)$
for $p\in S,a\in A$, and $F\in\Fil(A)$. The set $\Sblock$ of all
these carries a transition relation which induces discrete
furnishings $\JBB{\Sblock},\JDD{\Sblock}$ on $\Pcal\Sblock$. 
Next we convert $\Fcal\Agoth$ into a second garden $\Mcal\Agoth$, as
in the third column. This is essentially a furnished version of the
topology $\OS$. The carried furnishings $\JBB{\alpha},\JDD{\alpha}$
are lifted from the furnishings on $\Pcal\Sblock$ with
$$\proj^\la\circ\JBB{\alpha} \leq \JBB{\Sblock}\circ\proj^\la
\hspace{5ex}
\proj^\la\circ\JDD{\alpha} \leq \JDD{\Sblock}\circ\proj^\la$$
as the characterizing properties. The fourth column is
not needed until Subsection \ref{sec5.3}.

The end result of the 2-step process is the garden
$$\Mcal\Agoth = (\OS,\id,S)$$
which we think of as a manicured version of $\Agoth$. Throughout,  
the base space $S$ is unchanged. 

\TableAlg

Consider the unit $\Aunit_\Agoth$. Since everything happens
over the given base space $S$, it must be a frame morphism
$A\rTo\OS$ 
such that the triangle
\OBSCURE{
\Dm[size=2em]
A&            &\rTo&           &\OS   \\
 &\rdTo_\alpha&    &\ldTo_{\id}& \\
 &            &\OS &
\mD
}
\begin{center}
\includegraphics[scale=0.73,trim=0pt 0pt 0pt 0pt,clip]{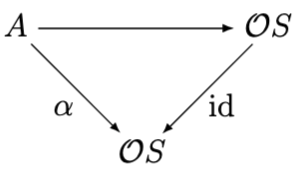}
\end{center}
commutes. Thus if there is such a
morphism then it can only be $\alpha$, the covering morphism of
$\Agoth$. This is certainly a frame morphism but we need a bit more. 

\Lem\label{item-250-a}
For each garden $\Agoth=(A,\alpha,S)$, the frame morphism $\alpha$ 
is a bed morphism relative to the harvest induced operations on $\OS$. 
\ELem

\Proof
Using the notation from above, we need to show that the comparisons
$$\alpha\circ\JBB{A} \leq \JBB{\alpha}\circ\alpha 
\hspace{7ex} 
  \alpha\circ\JDD{A} \leq \JDD{\alpha}\circ\alpha$$
hold. Thus, using the construction of $\JBB{\alpha},\JDD{\alpha}$ from
$\JBB{\Sblock},\JDD{\Sblock}$, it suffices to show that
$$\proj^\la\circ\alpha\circ\JBB{A}\leq\JBB{\Sblock}\circ\proj^\la\circ\alpha
\hspace{7ex}
\proj^\la\circ\alpha\circ\JDD{A}\leq\JDD{\Sblock}\circ\proj^\la\circ\alpha$$
hold. We look at these partly in parallel.

Consider any $x\in A$ and any flower
$$(p,a,F)\in\big(\proj^\la\circ\alpha\circ\JBB{A}\big)(x)
\hspace{7ex}
  (p,a,F)\in\big(\proj^\la\circ\alpha\circ\JDD{A}\big)(x)$$
that is a flower satisfying the two equivalent conditions
$$p\in\big(\alpha\circ\JBB{A}\big)(x) \hspace{3ex} 
   x\in\PBB{A}(p) 
\hspace{10ex} 
   p\in\big(\alpha\circ\JDD{A}\big)(x)  \hspace{3ex}
   x\notin\PDD{A}(p)$$
respectively. In both case we are concerned with transitions
$$(p,a,F)\rra(q,b,G)$$
in $\Sblock$. We must show $q\in\alpha x$, that is $x\in\nabla(q)$, for 
$$\mbox{all such transitions \hspace{5ex} at least one such transition}$$
respectively. 

For the left hand case we have a flower and a transition so that
$\PBB{A}(p) \subseteq F \subseteq \nabla(q)$ to give the required
result. 

For the right hand case since $(p,a,F)$ is a
flower we have $a\in\PDD{A}p$, so that $x\nleq a$, and hence the
$(-)$-healthiness of $\Sblock$ gives the required transition. \EP

This result justifies the following.

\Defn
For each garden $\Agoth=(A,\alpha,S)$ let
\OBSCURE{
\Dm \Agoth&\rTo^{\Aunit_\Agoth}&(\Gcal\circ\Fcal)\Agoth \mD
}
\begin{center}
\includegraphics[scale=0.73,trim=0pt 0pt 0pt 0pt,clip]{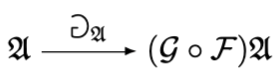}
\end{center}
be the garden map $(\id,\alpha)$. \ED
\EDefn

Now consider the naturality of $\Aunit_\Agoth$. Starting from a garden 
morphism
\OBSCURE{
\Dm \Bgoth&\rTo^{(\phi,f)}&\Agoth \mD
}
\begin{center}
\includegraphics[scale=0.73,trim=0pt 0pt 0pt 0pt,clip]{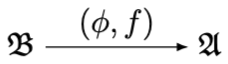}
\end{center}
between gardens $\Bgoth=(B,\beta,T)$ and $\Agoth=(A,\alpha,S)$, 
the functors $\Fcal$ and $\Gcal$ produce arrows
\OBSCURE{
\Dm[width=3.5em]
\Fcal\Agoth&\rTo^{(\phi,\Phi)}    &\Fcal\Bgoth& &
\Mcal\Bgoth&\rTo^{(\phi,\phi^\la)}&\Mcal\Bgoth
\mD
}
\begin{center}
\includegraphics[scale=0.73,trim=0pt 0pt 0pt 0pt,clip]{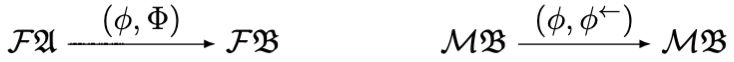}
\end{center}
where 
$$\Phi(p,a,F)=(\phi(p),f_\ast(a),f^\la(F))$$
for $(p,a,F)\in\Sblock$. We must show that the left hand square 
commutes. This expands 
\OBSCURE{
\Dm[size=2em]
               &                     &           &\hspace{7ex}& 
B              &                     &           &\rTo^\beta  &      
               &                     &\OT                     \\
               &                     &           &            &       
               &\rdTo~{\beta}        &           &            &      
               &\ldTo~{\id}                                   \\ 
\Bgoth         &\rTo^{(\beta,\id_T)} &\Mcal\Bgoth&            &      
               &                     &\OT        &\rTo^{\id}  &\OT  \\
\dTo~{(\phi,f)}&          &\dTo~{(\phi,\phi^\la)}&            &
\dTo^f         &                  &\dTo^{\phi^\la}&           &
\dTo_{\phi^\la}&                     &\dTo_{\phi^\la}         \\
\Agoth         &\rTo_{(\alpha,\id_S)}&\Mcal\Agoth&            &      
               &                     &\OS        &\rTo_{\id}  & \OS \\
               &                     &           &            &      
               &\ruTo~{\alpha}       &           &            &   
               &\luTo~{\id}                                   \\
               &                     &           &            &
A              &                     &           &\rTo_\alpha &       
               &                     &\OS                     \\
\mD
}
\begin{center}
\includegraphics[scale=0.73,trim=0pt 0pt 0pt 0pt,clip]{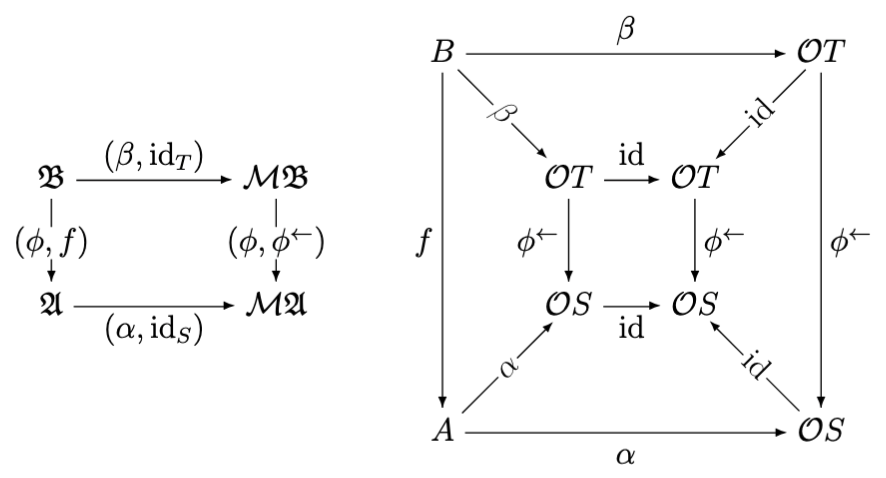}
\end{center}
to the right hand diagram, and of this, only showing the outer
square commutes is any problem. But this is precisely the property that
$(\phi,f)$ is a $\Grdn$-arrow. 

\Theom
The assignment $\Aunit_\bullet$ is a natural transformation 
\OBSCURE{
\Dm \Id&\rTo&\Gcal\circ\Fcal \mD
}
\begin{center}
\includegraphics[scale=0.73,trim=0pt 0pt 0pt 0pt,clip]{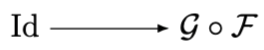}
\end{center}
from the identity functor $\Id$ to the endo-functor $\Gcal\circ\Fcal$
on $\Grdn$.  
\ETheom

\subsection{The geometric unit}\label{subsec5.2}

The geometric unit is not so straight forward.
For the endo-functor
$$\Hcal=\Fcal\circ\Gcal$$
on $\PlotL$, we attach to each plot $\Sgoth$ a lentile map 
\OBSCURE{
\Dm \Sgoth&\rTo^{\GunitS}&\Hcal\Sgoth \mD
}
\begin{center}
\includegraphics[scale=0.73,trim=0pt 0pt 0pt 0pt,clip]{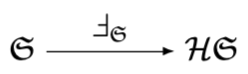}
\end{center}
To do this we fix some notation, as in Table
\ref{table-geo}, and review the relevant constructions.

\TableGeo

We start with a plot $\Sgoth$ as in the left-hand column of the table.
The transition relation on $\Sblock$ induces discrete
furnishings $\JBB{\Sblock},\JDD{\Sblock}$ on $\Pcal\Sblock$. The
functor $\Gcal$ converts this plot into a garden, as in the second
column. This is essentially a furnished version of the topology
$\OS$. The furnishings $\JBB{\sigma},\JDD{\sigma}$ are lifted from the
furnishings on $\Pcal\Sblock$ with
$$\sigma^\la\circ\JBB{\sigma} \leq \JBB{\Sblock}\circ\sigma^\la 
\hspace{5ex}
  \sigma^\la\circ\JDD{\sigma} \leq \JDD{\Sblock}\circ\sigma^\la$$
as the characterizing properties. 
Next we harvest $\Gcal\Sgoth$ to obtain
a second plot $\Hcal\Sgoth$, as in the third column. The transition
nodes, the healthy flowers, are certain triples 
$(p,U,\nabla)$
for $p\in S,U\in\OS$, and $\nabla\in\Fil(\OS)$
and these make up the new transition structure $\Sblockbar$.
To help with this
harvesting we attach a filter $\nabla(p)$ to each point $p\in S$. We
find that
$$U\in\nabla(p)\LRa p\in U \hspace{5ex}
 U\in\nabla(\sigma(P))\LRa P\in\sigma^\la(U)$$
for each $p\in S,P\in\Sblock$ and $U\in\OS$.
The transition relation on $\Sblockbar$ induces discrete
furnishings $\JBB{\Sblockbar},\JDD{\Sblockbar}$ on $\Pcal\Sblock$. 
The fourth column of the table is not needed until Subsection
\ref{sec5.3}.
Throughout this process the base
space $S$ is unchanged. 

Our first problem is to attach a flower to each parent transition
node $P\in\Sblock$.

\Defn
\label{item-250-g}
Let $\Sgoth$ be a plot, as above.
For each transition node $P\in\Sblock$ its {\sf root} is 
$\sigma(P)$, and its {\sf stalk} $\St(P)$ and {\sf bloom} $\Bl(P)$ are
given by 
$$U\subseteq \St(P) \LRa \sigma^\la(U)\subseteq P^\ra{}' \hspace{5ex} 
U\in\Bl(P) \LRa P^\ra\subseteq \sigma^\la(U)$$
(for $U\in\OS$). We set
$$\GunitS(P) = \big(\sigma(P), \St(P), \Bl(P)\big)$$
to produce a (potential) flower. \ED
\EDefn

Thus 
$\St(P) = (\sigma[P^\ra])^-{}'$
and $\Bl(P)$ is the filter of open neighbourhoods of
$\sigma[P^\ra]$.

\Lem
For each transition node $P$ (of the plot $\Sgoth$) the triple
$\GunitS(P)$ is a flower.
\ELem

\Proof
Let
$$\GunitS(P) = (p,U,\nabla)$$
where
$$p=\sigma(P) \hspace{3ex} U=\St(P) \hspace{3ex} \nabla = \Bl(P)$$
so we must show that
$$U\in\PDD{\OS}(p) \hspace{5ex} \PBB{\OS}(p)\subseteq\nabla$$
according to Definition \ref{item-240-b}.
We look at these separately.

Let $V=\JDD{\sigma}(U)$. By the construction of $\PDD{\OS}(p)$ we
have
$$U\in\PDD{\OS}(p) \LRa V\notin\nabla(p) \LRa p\notin V 
\LRa P\notin\sigma^\la(V)$$
so, by way of contradiction, suppose $P\in\sigma^\la(V)$. 
We have
$$P\in\sigma^\la(V) \subseteq 
\big(\JDD{\Sblock}\circ\sigma^\la\big)(U)
\hspace{5ex}
\sigma^\la(U) \subseteq P^\ra{}'$$
by the definitions of $\JDD{\sigma}$ and $\St P$. The left hand
inclusion gives some node $Q$ with
$$P\rra Q \in \sigma^\la(U) \subseteq P^\ra{}'$$
which is patent nonsense.

Remembering the definition of $\PBB{\OS}$, for each $W\in\OS$ we have
$$W\in\PBB{\OS}(p)
\LRa \JBB{\sigma}(W)\in\nabla(p)
\LRa p\in \JBB{\sigma}(W)
\LRa P\in(\sigma^\la\circ\JBB{\sigma})(W)$$
so that
$$W\in\PBB{\OS}(p)
\Rra P\in\big(\sigma^\la\circ\JBB{\sigma}\big)(W)
\subseteq\big(\JBB{\Sblock}\circ\sigma^\la\big)(W)
\Rra P^\ra\subseteq\sigma^\la(W)
\Rra W\in\Bl(P)$$
to give the required result. \EP

Our next job is to show that these flowers $\GunitS(P)$ are healthy. 

\Lem\label{item-250-f}
We have
$$P\rra Q \Rra \GunitS(P)\rra \GunitS(Q)$$
for transition nodes $P,Q$ (of the parent plot $\Sgoth$). 
\ELem

\Proof
Assuming $P\rra Q$ we must check that 
$$\St(P)\notin\big(\nabla\circ\sigma\big)(Q)
\hspace{5ex}
  \Bl(P)\subseteq\big(\nabla\circ\sigma\big)(Q)$$
hold. 
 
For the left hand condition, with 
$U=\St(P)$ we have $\sigma^\la(U)\subseteq
P^\ra{}'$, to give $Q\notin\sigma^\la(U)$, that is $\sigma(Q)\notin U$,
and hence $U\notin\big(\nabla\circ\sigma\big)(Q)$.

For the right hand condition, with $U\in\Bl(P)$ we have
$Q\in P^\ra\subseteq\sigma^\la(U)$, so that $\sigma(Q)\in U$, and hence
$U\in\big(\nabla\circ\sigma\big)(Q)$, as required.  \EP

With this we can give the flowers the required certificate of health.

\Lem
For each plot $\Sgoth$ the set of flowers of the form $\GunitS(P)$
(for $P\in\Sblock$) is healthy. 
\ELem

\Proof
$(\upa)$ Consider any transition node $P$ and open set 
$U\notin\Bl(P)$. Then $P^\ra\nsubseteq\sigma^\la(U)$, to give the
required witness.  

$(-)$ Consider any transition node $P$ and open set 
$U\nsubseteq\St(P)$. Then $\sigma^\la(U)\nsubseteq P^\ra{}'$ and hence
there is some $Q\in\sigma^\la(U)$ with $P\rra Q$. Lemma
\ref{item-250-f} gives $\GunitS(P)\rra\GunitS(Q)$ and this with
$\sigma(Q)\in U$ gives the required witness. \EP

Within the harvest $\Hcal\Sgoth$ of the
garden $\Gcal\Sgoth$ obtained from a plot $\Sgoth$ we have located a
set of healthy flowers, namely those arising from the parent transition nodes. 

\Theom
For each plot $\Sgoth$ the assignment
\OBSCURE{
\Dm[height=1.3em]
\Sgoth&\rTo&\Hcal\Sgoth \\
 P    &\rMapsto& \GunitS(P)
\mD
}
\begin{center}
\includegraphics[scale=0.73,trim=0pt 0pt 0pt 0pt,clip]{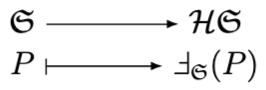}
\end{center}
is a lentile map.
\ETheom

\Proof
By Lemma \ref{item-250-f} this assignment is a plot map.
Thus it suffices to check the properties $(\upa,-)$ of Lemma
\ref{item-220-g}. For each we start from a node $P$ of $\Sgoth$ and a
transition
$$\big(\sigma(P),\St(P),\Bl(P)) \rra (r,W,H)$$
in the healthy harvest, and we must produce a certain transition
$P\rra Q$ in $\Sgoth$.

$(\upa)$ The given transition ensures that 
$\Bl(P)\subseteq\nabla(r)$, and hence
$$P^\ra\subseteq \sigma^\la(U) \Rra U\in\Bl(P) \Rra r\in U$$
for each $U\in\OS$. In particular, we have 
$P^\ra\subsetneq \sigma^\la(r^-{}')$ which gives a transition 
$P\rra Q$ with $\sigma(Q)\in r^-$, as required.

$(-)$ As well as the transition we are given some $r\in U\in \OS$. The
transition ensures that  
$r\notin\St(P)$, so that $U\subsetneq\St(P)$, and hence $\sigma^\la(U)$
meets $P^\ra$. This gives a transition $P\rra Q$ with $\sigma(Q)\in U$,
as required. \EP

This completes the first objective.  The second objective is to show
that $\Gunit_\bullet$ is a natural transformation from the identity
functor to the endo-functor $\Fcal\circ\Gcal$ on $\PlotL$. In more
detail we show that for each lentile map $(\phi,\Phi)$ the right hand
square commutes.
\OBSCURE{
\Dm[size=2em] & & &\hspace{10ex}&
\Sgoth            &\rTo^{\GunitS}    &\Hcal\Sgoth             \\
\Sgoth &\rTo^{(\phi,\Phi)}&\Tgoth &&
\dTo^{(\phi,\Phi)}&                  &\dTo_{\Hcal(\phi,\Phi)} \\
& & && \Tgoth &\rTo_{\Gunit_\Tgoth} &\Hcal\Tgoth \mD
}
\begin{center}
\includegraphics[scale=0.73,trim=0pt 0pt 0pt 0pt,clip]{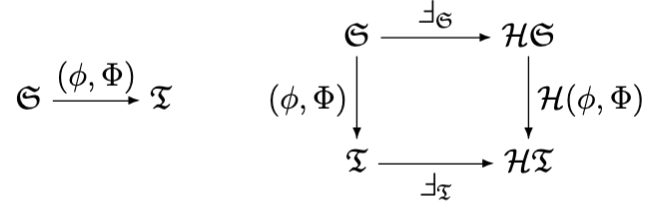}
\end{center}

We begin by taking a more detailed look at the components of the
putative commuting square. This will help to set up the notation
used in the analysis.

Each pair of plots
$$\Sgoth=(\Sblock,\sigma,S) \hspace{5ex} \Tgoth=(\Tblock,\tau,T)$$
gives a pair of gardens
$$\Gcal\Sgoth=(\OS,\id,S)\hspace{5ex}\Gcal\Tgoth=(\OT,\id,T)$$
each of which is a furnished version of the topology on the
corresponding base space. The source $\PlotL$-arrow $(\phi,\Phi)$ gives
a garden morphism
\OBSCURE{
\Dm[width=17ex]
\Gcal\Sgoth&\lTo^{\Gcal(\phi,\Phi)=(\phi,f=\phi^\la)}&\Gcal\Tgoth
\mD
}
\begin{center}
\includegraphics[scale=0.73,trim=0pt 0pt 0pt 0pt,clip]{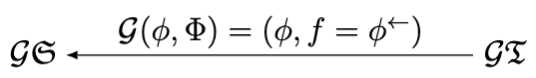}
\end{center}
where the frame morphism
\OBSCURE{
\Dm[width=8ex] \OS&\lTo^{f=\phi^\la}&\OT \mD
}
\begin{center}
\includegraphics[scale=0.73,trim=0pt 0pt 0pt 0pt,clip]{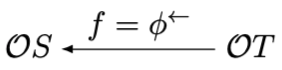}
\end{center}
is the inverse image function of the continuous map
$\phi$ between the base spaces. 
We now harvest each of the two gardens to produce a new pair
$$\Hcal\Sgoth=(\Fcal\circ\Gcal)\Sgoth=(\Sblock,\proj,S) \hspace{3ex}
  \Hcal\Tgoth=(\Fcal\circ\Gcal)\Tgoth=(\Tblock,\proj,T)$$
of plots. The corresponding transition nodes, the members of
$\Sblock$ and $\Tblock$, are certain triples
$$(p,U,F) \hspace{7ex} (q,V,G)$$
where
$$p\in S \hspace{2ex} U\in\OS \hspace{2ex} F\in\Fil(\OS)
\hspace{5ex}
  q\in T \hspace{2ex} V\in\OT \hspace{2ex} G\in\Fil(\OT)$$
are the components. These are subject to certain
restrictions. Each source node
$$P\in\Sblock \hspace{5ex} Q\in\Tblock$$
gives a target node 
$$\GunitS(P)       =\big(\sigma(P),\St(P),\Bl(P)\big) 
\hspace{5ex}
  \Gunit_\Tgoth(Q) =\big(\tau(Q),\St(Q),\Bl(Q)\big)$$
in the corresponding harvest plot. 
The garden morphism $\Gcal(\phi,\Phi)$ gives a $\PlotL$-arrow
\OBSCURE{
\Dm[width=15ex]
\Sblockbar &\rTo^{(\phi,\, \phi_\ast=f_\ast,\,\phi^\La=f^\la)}
                   &\overline{\Tblock}
\mD
}
\begin{center}
\includegraphics[scale=0.73,trim=0pt 0pt 0pt 0pt,clip]{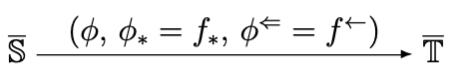}
\end{center}
over the given map $\phi$. We look at
the three components of this, and explain the
notation. 

The root component $\phi$ is the parent continuous map.

The stalk component $\phi_\ast=f_\ast$ is the right adjoint
\OBSCURE{
\Dm[width=8ex]
\OS&\pile{\lTo^{f=\phi^\la} \\ \\ \rTo_{\phi_\ast=f_\ast}}&\OT
\mD
}
\begin{center}
\includegraphics[scale=0.73,trim=0pt 0pt 0pt 0pt,clip]{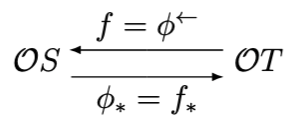}
\end{center}
of the inverse image map. In particular, we have
$$V\subseteq\phi_\ast(U) \LRa \phi^\la(V)\subseteq U$$
for $U\in\OS$ and $V\in\OT$.

The bloom component
\OBSCURE{
\Dm[width=8ex] \Fil(\OS)&\rTo^{\phi^\La=f^\la}&\Fil(\OT) \mD
}
\begin{center}
\includegraphics[scale=0.73,trim=0pt 0pt 0pt 0pt,clip]{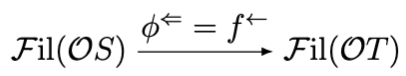}
\end{center}
is the inverse image map of $f=\phi^\la$ which itself is an inverse
image map. Thus we have
$$V\in\phi^\La(F) \LRa  \phi^\la(V)\in F$$
for each $F\in\Fil(\OS)$ and $V\in\OT$.

\TableNat

By unravelling the square we obtain the diagram of functions on the
left of Table \ref{table-205}. This must commute. 
Only the outer square is
problematic, so we concentrate on that.  
By tracking a node $P\in\Sblock$ around the square, as on the right of
Table \ref{table-205},  
we see that
$$\mbox{(r)}\hspace{2ex}\tau\circ\Phi=\phi\circ\sigma  \hspace{5ex}
  \mbox{(s)}\hspace{2ex} \St\circ\Phi=\phi_\ast\circ\St \hspace{5ex}
  \mbox{(b)}\hspace{2ex} \Bl\circ\Phi=\phi^\La\circ\Bl$$
are the three required identities. The first of these, (r), is given
(since $(\phi,\Phi)$ is a plot map), so it remains to deal with
(s) and (b). 

\Theom
The assignment $\Gunit_\bullet$ is a natural transformation
\OBSCURE{
\Dm \Id\Tlens&\rTo&\Fcal\circ\Gcal \mD
}
\begin{center}
\includegraphics[scale=0.73,trim=0pt 0pt 0pt 0pt,clip]{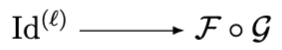}
\end{center}
from the identity functor $\Id\Tlens$ to the endo-functor
$\Fcal\circ\Gcal$ on $\PlotL$.  
\ETheom

\Proof
As indicated above, it suffices to verify the two
equalities (s,b). 

\smallskip

(s) For each $P\in\Sblock$ and $V\in\OT$ we have
$$V\subseteq(\St\circ\Phi)(P) \LRa \tau^\la(V)\subseteq(\Phi(P))^\ra{}'$$
and
$$V\subseteq(\phi_\ast\circ\St)(P) 
\LRa \phi^\la(V)\subseteq \St(P)
\LRa (\sigma^\la\circ\phi^\la)(V)\subseteq P^\ra{}'
\LRa (\Phi^\la\circ\tau^\la)(V)  \subseteq P^\ra{}'$$
so that an equivalence
$$\tau^\la(V)\mbox{ meets }(\Phi(P))^\ra \LRa
(\Phi^\la\circ\tau^\la)(V)\mbox{ meets }P^\ra$$
is the essential content of the required equality.

We deal separately with the two directions of the equivalence.

Suppose $\tau^\la(V)$ meets $(\Phi(P))^\ra$, 
so that
$$\Phi P\rra R \hspace{5ex} \tau R\in V$$
for some node $R\in\Tblock$. Since $(\phi,\Phi)$ is a lentile map 
we have
$$P\rra Q \hspace{5ex} (\tau\circ\Phi)(Q)=(\phi\circ\sigma)(Q)\in V$$
for some node $Q\in \Sblock$,
so that
$(\Phi^\la\circ\tau^\la)(V)$ meets $P^\ra$ 
at $Q$.

Conversely, suppose
$(\Phi^\la\circ\tau^\la)(V)$ meets $P^\ra$ 
so that
$$P\rra Q \hspace{5ex} (\tau\circ\Phi)(Q)\in V$$
for some node $Q\in\Sblock$. Since $\Phi$ is a transition morphism
this gives $\Phi P\rra \Phi Q$ and hence
$\tau^\la(V)$  meets $(\Phi P)^\ra$ at $\Phi Q$.

\smallskip

(b) We show 
$$(\Bl\circ\Phi)(P) = (\phi^\La\circ\Bl)(P)$$
for arbitrary $P\in\Sblock$. For each $V\in\OT$ we have
$$V\in(\Bl\circ\Phi)(P) \LRa (\Phi(P))^\ra\subseteq\tau^\la(V)$$
and
$$V\in(\phi^\La\circ\Bl)(P)
\LRa \phi^\la(V)\in\Bl(P)
\LRa P^\ra\subseteq (\sigma^\la\circ\phi^\la)(V)=(\Phi^\la\circ\tau^\la)(V)$$
so that an equivalence
$$(\Phi(P))^\ra\subseteq\tau^\la(V)
\LRa
P^\ra\subseteq (\sigma^\la\circ\phi^\la)(V) = (\Phi^\la\circ\tau^\la)(V)$$
is the essential content of the required equality.

Suppose $(\Phi P)^\ra\subseteq\tau^\la(V)$ and consider a
transition $P\rra Q$ in $\Sblock$. 
This gives a transition $\Phi(P)\rra\Phi(Q)$ in $\Tblock$, so that
$\Phi(Q)\in\tau^\la(V)$, 
as required.

Conversely, suppose
$P^\ra\subseteq (\sigma^\la\circ\phi^\la)(V)$
and consider any transition $\Phi(P)\rra R$ in $\Tblock$. Since we are
dealing with a lentile map this gives 
$$P\rra Q \hspace{5ex} (\phi\circ\sigma)(Q)\sqsubseteq\tau(R)$$
for some $Q\in\Sblock$. But now $Q\in(\sigma^\la\circ\phi^\la)(V)$ so
that $\tau(R)\in V$, as required. \ED

\subsection{Idempotency}\label{sec5.3}

The results so far do not show that we have an
adjunction. We certainly have a pair of functors $\Fcal$ and $\Gcal$,
and a pair of natural transformations
\OBSCURE{
\Dm 
\Id_{\Grdn}&\rTo^{\Aunit}&\Gcal\circ\Fcal &&
\Id_{\PlotL}&\rTo^{\Gunit}&\Fcal\circ\Gcal 
\mD
}
\begin{center}
\includegraphics[scale=0.73,trim=0pt 0pt 0pt 0pt,clip]{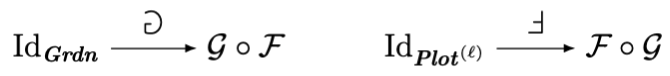}
\end{center}
where the two source endo-functors are the
appropriate identity functors. For an adjunction we require certain
triangle laws to hold, but we show more.

Starting from a garden $\Agoth$ and a plot $\Sgoth$ 
we may hit each unit $\Aunit_\Agoth$ and $\GunitS$
with the appropriate functor
to transpose it to the other side
\OBSCURE{
\Dm 
\Gcalbar\Sgoth&\rTo^{\Gcal(\GunitS)}   &\Gcal\Sgoth  &&
\Fcalbar\Agoth&\rTo^{\Fcal(\Aunit_\Agoth)}&\Fcal\Agoth
\mD
}
\begin{center}
\includegraphics[scale=0.73,trim=0pt 0pt 0pt 0pt,clip]{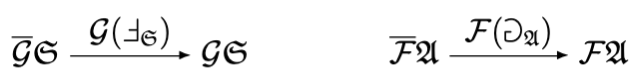}
\end{center}
where
$$\Gcalbar= \Gcal\circ\Fcal\circ\Gcal \hspace{9ex}
  \Fcalbar= \Fcal\circ\Gcal\circ\Fcal$$
are the two source functors. As particular cases of the units
we have
\OBSCURE{
\Dm 
\Gcal\Sgoth&\rTo^{\Aunit_{\Gcal\Sgoth}}&\Gcalbar\Sgoth &&
\Fcal\Agoth&\rTo^{\Gunit_{\Fcal\Agoth}}&\Fcalbar\Agoth  
\mD
}
\begin{center}
\includegraphics[scale=0.73,trim=0pt 0pt 0pt 0pt,clip]{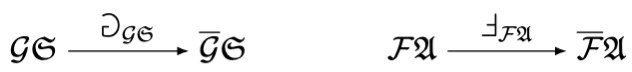}
\end{center}
and an adjunction requires that each of the two composites
\OBSCURE{
\Dm 
\Gcal\Sgoth&\rTo^{\Aunit_{\Gcal\Sgoth}}&\Gcalbar\Sgoth 
           &\rTo^{\Gcal(\GunitS)}   &\Gcal\Sgoth    &&
\Fcal\Agoth&\rTo^{\Gunit_{\Fcal\Agoth}}   &\Fcalbar\Agoth  
           &\rTo^{\Fcal(\Aunit_\Agoth)}&\Fcal\Agoth
\mD
}
\begin{center}
\includegraphics[scale=0.73,trim=0pt 0pt 0pt 0pt,clip]{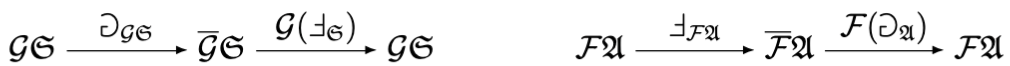}
\end{center}
is the identity on its end object. We show that all four possible 
composites are identities, 
so that each transposed unit is an isomorphism, and we have an idempotent
adjunction.

\subsubsection*{The algebraic induced isomorphism}

We continue with the notation set out in Table \ref{table-alg}. The
first three columns are explained in Subsection \ref{subsec5.1}. We now also
use the fourth column. Thus we harvest the garden 
$\Mcal\Agoth$ to obtain a second plot
$\Fcalbar\Agoth$, as in the right-hand column. 
The transition nodes are certain triples
$(p,U,\nabla)$
for $p\in S,U\in\OS$, and $\nabla\in\Fil(\OS)$. The set
$\Sblockbar$ of all of these carries a transition
relation. We don't need the induced furnishings on $\Pcal\Sblockbar$.

The heart of our problem is to relate the
nodes
$$(p,a,F)\in\Sblock \hspace{5ex} (p,U,\nabla)\in\Sblockbar$$
of the two constructed plots. To do that we remember a little bit of
basic frame theory. 

The covering morphism $\alpha^\ast=\alpha$ is a frame morphism
and, as such, has a right adjoint
\OBSCURE{
\Dm A&\pile{\rTo^{\alpha^\ast} \\ \\ \lTo_{\alpha_\ast}}&\OS \mD
}
\begin{center}
\includegraphics[scale=0.73,trim=0pt 0pt 0pt 0pt,clip]{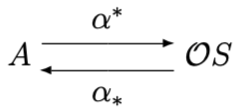}
\end{center}
where 
the adjunction property is
$$\alpha^\ast(x) \subseteq U \LRa x\leq \alpha_\ast(U)$$
for $x\in A$ and $U\in\OS$.
The composite $j_\Agoth=\alpha_\ast\circ\alpha^\ast$ is a closure operation
on $A$. The 
surjectivity of $\alpha$ ensures the other composite
$\alpha^\ast\circ\alpha_\ast$ is the identity on $\OS$. Thus we have
$$x\leq j_\Agoth(x) \hspace{5ex} (\alpha^\ast\circ\alpha_\ast)(U)=U$$
for $x\in A$ and $U\in\OS$.

There is a similar connection between the filters on $A$ and $\OS$
\OBSCURE{
\Dm[height=1.5em,width=4em]
   F    &\rMapsto           & \alpha[F] \\
\Fil(A) &\pile{\rTo \\ \\ \lTo}& \Fil(\OS) \\
\alpha^\la(\nabla)&\lMapsto   & \nabla
\mD
}
\begin{center}
\includegraphics[scale=0.73,trim=0pt 0pt 0pt 0pt,clip]{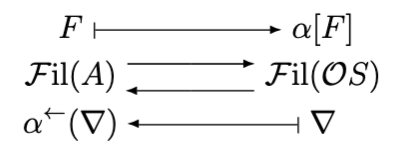}
\end{center}
obtained by taking direct and inverse images, as indicated. The
composite $J_\Agoth=\alpha^\la\circ\alpha[\cdot]$ inflates each filter
$F\in\Fil(A)$, and the other composite is the identity.

Each flower $P=(p,a,F)$ of the harvest of $\Agoth$ has a stalk 
$a\in A$ and a bloom $F\in\Fil(A)$. These can not be chosen
arbitrarily. 

\Lem\label{item-250-j}
For each flower  $P=(p,a,F)\in\Sblock$ we have 
$j_\Agoth(a) = a$ and $J_\Agoth(F)=F$.

\ELem

\Proof
Let $x=j_\Agoth(a)$, so that $\alpha(x)=\alpha(a)$ and $a\leq x$. We must
show that $x\leq a$. 

If $x\nleq a$ then, since $(p,a,F)$ is healthy, there is a transition 
$$(p,a,F)\rra(q,b,G)$$
in $\Sblock$ with $x\in\nabla(q)$, that is $q\in\alpha(x)$. But this
transition requires $a\notin\nabla(q)$, that is $q\notin\alpha(a)$, which
is contradictory since $\alpha(x)=\alpha(a)$. 

On general grounds we have $F\subseteq J_\Agoth(F)$, so we require the
converse inclusion.
Consider any $x\in J_\Agoth(F)$. We have $\alpha(x)\in\alpha[F]$, so
that $\alpha(x)=\alpha(y)$ for some $y\in F$. By way of contradiction
suppose $x\notin F$. The parent node $P=(p,a,F)$ is healthy,
so there is a second transition, as above, this time 
with $x\notin\nabla(q)$, that is $q\notin\alpha(x)$. This
transition requires $F\subseteq\nabla(q)$, so that $y\in\nabla(q)$,
that is $q\in\alpha(y)$. This contradicts $\alpha(x)=\alpha(y)$. \EP

To show that $\Gunit_{\Fcal\Agoth}$ and $\Fcal(\Aunit_\Agoth)$ form an
inverse pair we describe the behaviour of each of these assignments,
beginning with $\Gunit_{\Fcal\Agoth}$.

Each transition node $P$ of
$\Fcal\Agoth$, that is each $P\in\Sblock$, has the form
$$P=(p,a,F)$$
for some $p\in S,a\in A,F\in\Fil(A)$. Then
$$\Gunit_{\Fcal\Agoth}P = (\overline{p},U,\nabla)$$
where, by Definition \ref{item-250-g} we have
$$\overline{p}=\proj (P) = p\hspace{5ex}U=\St(P)\hspace{5ex}\nabla=\Bl(P)$$
and the stalk and bloom are given by
$$V\subseteq \St(P) \LRa \proj^\la(V)\subseteq P^\ra{}' \hspace{5ex}
  V\in\Bl(P) \LRa P^\ra \subseteq \proj^\la(V)$$
for each $V\in\OS$. The trick we use is that each such $V$ has the
form $\alpha^\ast(x)$ for some $x\in A$ (for instance 
$x=\alpha_\ast(V)$). There is also another crucial observation.

\Theom\label{item-250-k}
For each garden $\Agoth$ we have 
$$\Gunit_{\Fcal\Agoth}P=(p,\alpha(a),\alpha[F])$$
for each transition node $P=(p,a,F)$ of $\Fcal\Agoth$.
\ETheom

\Proof
We make use of the notation above.

\smallskip

(Stalk) Let $U=\St(P)$. We first show that 
$$x\leq a \LRa \alpha(x)\subseteq U$$
holds for all $x\in A$. In fact, we show the contrapositive.

Suppose $\alpha(x)\nsubseteq U$. Then, by the construction of $\St(P)$, we
see that $(\proj^\la\circ\alpha)(x)$ meets $P^\ra$ to give some transition
$$(p,a,F)\rra(q,b,G)\in(\proj^\la\circ\alpha)(x)$$
in $\Sblock$. This requires
$$q\in\alpha(x) \hspace{5ex} q\notin\alpha(a)$$
so that $\alpha(x)\nsubseteq\alpha(a)$, and hence $x\nleq a$.

Conversely, suppose $x\nleq a$. Since $P$ is healthy there is a
transition
$$(p,a,F)\rra(q,b,G)$$
in $\Sblock$ with $x\in\nabla(q)$. But then $\alpha(x)\nsubseteq U$ by
the converse of the previous argument. 

To complete this part of the proof consider any $V\in\OS$. Since
$\alpha$ is surjective we 
have $V=\alpha(x)$ for some $x\in A$. With this we have
$$V\subseteq\St(P) \Rra x\leq a \Rra \alpha(x)\subseteq \alpha(a) \Rra 
V\subseteq \alpha(a)$$
using above equivalence at the first step. Conversely
$$V\subseteq\alpha(a) \Rra \alpha(x)\subseteq \alpha(a) \Rra 
x\leq j(a)=a \Rra V=\alpha(x)\subseteq\St(P)$$
using Lemma \ref{item-250-j} at the second step and the above equivalence at
the third. 

Thus we have
$$V\subseteq\St(P)\LRa V\subseteq \alpha(a)$$
to give $\St (P) = \alpha (a)$, as required. 

\smallskip

(Bloom) Let $\nabla=\Bl(P)$. Since $\alpha$ is surjective we see that
an equivalence
$$x\in F \LRa \alpha(x)\in \nabla$$
(for $x\in A$) will lead to the required result. As before, we prove
the contrapositive.

Suppose $\alpha (x) \notin\nabla$. 
Then, by the construction of $\Bl(P)$, we see that 
$$P^\ra \nsubseteq (\proj^\la\circ\alpha)(x)$$
to give some transition
$$(p,a,F)\rra(q,b,G)\notin(\proj^\la\circ\alpha)(x)$$
in $\Sblock$. This requires
$$q\notin \alpha(x) \hspace{5ex} F\subseteq\nabla(q)$$
so that $x\notin F$, as required. 

Conversely, suppose $x\notin F$. Since $P$ is healthy there is a
transition
$$(p,a,F)\rra(q,b,G)$$
in $\Sblock$ with $x\notin\nabla(q)$. This transition witnesses that
$$P^\ra \nsubseteq (\proj^\la\circ\alpha)(x)$$
and hence $\alpha (x) \notin\nabla$, to complete the whole proof. \EP

To get a description of the map
$\Fcal(\Aunit_\Agoth)$ we remember what the algebraic unit
\OBSCURE{
\Dm \Agoth&\rTo^{\Aunit_\Agoth}&\Mcal\Agoth \mD
}
\begin{center}
\includegraphics[scale=0.73,trim=0pt 0pt 0pt 0pt,clip]{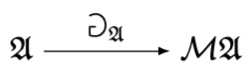}
\end{center}
look likes. By Subsection \ref{subsec5.1} this is essentially the covering
morphism $\alpha$ of $\Agoth$ 
viewed as an enriched morphism. Thus, by Subsection \ref{subsec4.2} the map
$\Fcal(\AunitA)$ is given by
\OBSCURE{
\Dm[height=1.25em,width=15ex]
\Fcalbar\Agoth&\rTo&\Fcal\Agoth \\
(p,U,\nabla)&\rMapsto&(p,\alpha_\ast(U),\alpha^\la(\nabla))
\mD
}
\begin{center}
\includegraphics[scale=0.73,trim=0pt 0pt 0pt 0pt,clip]{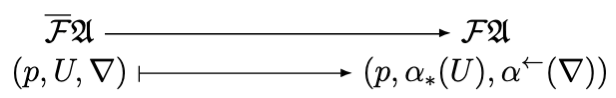}
\end{center}
for $p\in S,U\in\OS,\nabla\in\Fil(\OS)$ where the source triple
must satisfy certain conditions.
Since $\alpha$ is surjective we can always represent $U$ as 
$\alpha (a)$ and $\nabla$ as $\alpha[F]$. Thus
\OBSCURE{
\Dm[height=1.25em,width=15ex]
\Fcalbar\Agoth&\rTo&\Fcal\Agoth \\
(p,\alpha(a),\alpha[F])&\rMapsto&(p,j_\Agoth(a),J_\Agoth(F))
\mD
}
\begin{center}
\includegraphics[scale=0.73,trim=0pt 0pt 0pt 0pt,clip]{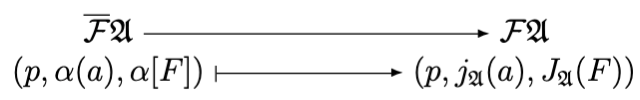}
\end{center}
is a description of the transpose of this unit. 
Even better, we can always choose the representatives $a,F$ with
$j_\Agoth(a)= a$ and $J_\Agoth(F)=F$, so the
transposed unit is
\OBSCURE{
\Dm[height=1.25em,width=15ex]
\Fcalbar\Agoth&\rTo&\Fcal\Agoth \\
(p,\alpha(a),\alpha[F])&\rMapsto&(p,a,F)
\mD
}
\begin{center}
\includegraphics[scale=0.73,trim=0pt 0pt 0pt 0pt,clip]{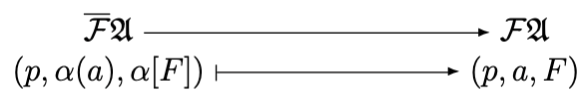}
\end{center}
and we are beginning to see what is going on. 

\Theom
For each garden $\Agoth$ the two assignments
\OBSCURE{
\Dm 
   \Fcal\Agoth&\rTo^{\Gunit_{\Fcal\Agoth}}&\Fcalbar\Agoth&&
\Fcalbar\Agoth&\rTo^{\Fcal(\Aunit_\Agoth)}&\Fcal\Agoth
\mD
}
\begin{center}
\includegraphics[scale=0.73,trim=0pt 0pt 0pt 0pt,clip]{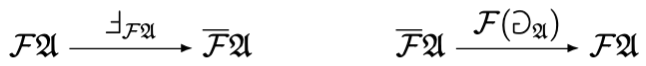}
\end{center}
form an inverse pair is isomorphisms.
\ETheom

\Proof
We continue with the notation developed above.

On general grounds both the assignments are the transition part of a
lentile map. Thus it suffices that the two composite are the
appropriate identity function.  

By Theorem \ref{item-250-k} and the remarks above the 
composite via $\Fcalbar\Agoth$ is
\OBSCURE{
\Dm[height=1.25em,width=10ex]
\Fcal\Agoth&\rTo&\Fcalbar\Agoth&\rTo&\Fcal\Agoth \\
(p,a,F)&\rMapsto&(p,\alpha(a),\alpha[F])&\rMapsto&(p,j_\Agoth(a),J_\Agoth(F))
\mD
}
\begin{center}
\includegraphics[scale=0.73,trim=0pt 0pt 0pt 0pt,clip]{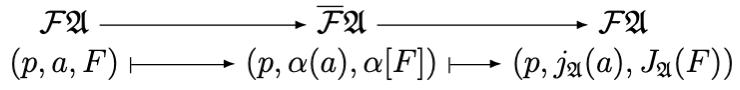}
\end{center}
for $P=(p,a,F)\in\Sblock$. Lemma \ref{item-250-j} shows that this
composite is the identity on $\Sblock$.

\smallskip

For the composite via $\Fcal\Agoth$, consider any node
$P=(p,U,\nabla)\in\Sblockbar$. As remarked above we have $U=\alpha (a)$
for some $a\in A$ with $j_\Agoth(a)=a$, and $\nabla=\alpha[F]$ for some
$F\in\Fil(A)$ with $J_\Agoth(F)=F$. Thus, using Theorem \ref{item-250-k} 
we see that the composite assignment is
\OBSCURE{
\Dm[height=1.25em,width=10ex]
\Fcalbar\Agoth&\rTo&\Fcal\Agoth&\rTo&\Fcalbar\Agoth \\
(p,\alpha(a),\alpha[F])&\rMapsto&(p,a,F)&\rMapsto&(p,\alpha(a),\alpha[F])
\mD
}
\begin{center}
\includegraphics[scale=0.73,trim=0pt 0pt 0pt 0pt,clip]{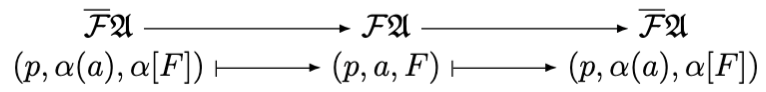}
\end{center}
which gives the required result. \ED

\subsubsection*{The geometric induced isomorphism}

For each plot $\Sgoth$ the constructions give two
garden morphisms
\OBSCURE{
\Dm 
   \Gcal\Sgoth&\rTo^{\Aunit_{\Gcal\Sgoth}}&\Gcalbar\Sgoth&&
\Gcalbar\Sgoth&\rTo^{\Gcal(\GunitS)}      &\Gcal\Sgoth
\mD
}
\begin{center}
\includegraphics[scale=0.73,trim=0pt 0pt 0pt 0pt,clip]{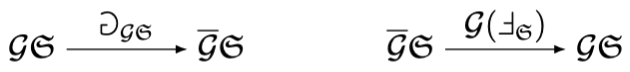}
\end{center}
in opposite directions. To form an adjunction the composite 
through $\Gcalbar\Sgoth$ must be the
identity on $\Gcal\Sgoth$. In this block we show that both composites
are the identity, and hence the two maps form an inverse pair of isomorphisms.
In fact, we show more. Each of the two components is an identity, that
is the two gardens $\Gcal\Sgoth$ and $\Gcalbar\Sgoth$ are exactly the
same. Curiously, the proof of this is somewhat easier than the proof
of the corresponding result for gardens given in the previous block. 

We continue with the notation of Table \ref{table-geo}. The
first three columns are explained earlier. We now also
use the fourth column. The transition relation on $\Sblockbar$ induces discrete
furnishings $\JBB{\Sblockbar},\JDD{\Sblockbar}$ on $\Pcal\Sblockbar$. 
We lift these to obtain another furnished version of $\OS$.

For these constructions we have a fixed base space
$S$, so the topological component of each arrow we meet is the
identity on $S$. It's the other components we need to sort out. 

For any garden $\Agoth$ the unit
\OBSCURE{
\Dm \Agoth&\rTo^{\Aunit_\Agoth}&\Mcal\Agoth \mD
}
\begin{center}
\includegraphics[scale=0.73,trim=0pt 0pt 0pt 0pt,clip]{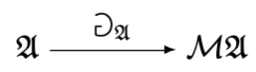}
\end{center}
is essentially the covering morphism of $\Agoth$. By Lemma \ref{item-250-a}
this is always an enriched morphism from the parent furnished frame to the
furnished version of the topology of the base space. For the particular case
$\Agoth=\Gcal\Sgoth$ the covering morphism is the identity on
$\OS$. Thus we obtain the following.

\Lem
For each plot $\Sgoth$ we have
$$\JBB{\sigma}\leq\JBB{\sigmaB} \hspace{5ex}
\JDD{\sigma}\leq\JDD{\sigmaB}$$
(in the notation above).
\ELem

This result deals with $\Aunit_{\Gcal\Sgoth}$. For $\Gcal(\GunitS)$
we remember that a use of $\Gcal$ always produces an enriched morphism
between the two furnished topologies. For the unit
\OBSCURE{
\Dm \Sgoth&\rTo^\GunitS&\Hcal\Sgoth \mD
}
\begin{center}
\includegraphics[scale=0.73,trim=0pt 0pt 0pt 0pt,clip]{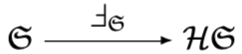}
\end{center}
the two topologies are the same with the identity morphism passing
between them. The two sets of furnishings may be different, but the
behaviour of $\Gcal$ gives the following.

\Lem
For each plot $\Sgoth$ we have
$$\JBB{\sigmaB}\leq\JBB{\sigma} \hspace{5ex}
\JDD{\sigmaB}\leq\JDD{\sigma}$$
(in the notation above).
\ELem

These two lemmas combine to give the main result of this block.

\Theom
For each plot $\Sgoth$ the two gardens $\Gcal\Sgoth$ and
$\Gcalbar\Sgoth$ are exactly the same.
\ETheom

Thus we have set up an idempotent contravariant adjunction.

\bibliographystyle{plain}
\bibliography{topsyrefs}
\end{document}